\documentclass[11pt]{amsart}
\usepackage{amsmath}
\usepackage{amsthm}
\usepackage{amsbsy}
\usepackage{amsopn}
\usepackage{amsmath,amscd}
\usepackage{amssymb}
\usepackage{amsfonts}
\usepackage{multirow}
\usepackage{mathrsfs} 
\usepackage[official ]{eurosym}
\usepackage[backref=true, maxbibnames=4,style=alphabetic, url=false, isbn=false, doi=false]{biblatex}

\newcommand{\RNum}[1]{\uppercase\expandafter{\romannumeral #1\relax}}

\newcommand\cA{\mathcal{A}}

\newcommand\cC{\mathcal{C}}
\newcommand\cE{\mathcal{E}}
\newcommand\cF{\mathcal{F}}
\newcommand\cG{\mathcal{G}}
\newcommand\cH{\mathcal{H}}
\newcommand\cI{\mathcal{I}}

\newcommand\cK{\mathcal{K}}

\newcommand\cO{\mathcal{O}}
\newcommand\cP{\mathcal{P}}

\newcommand\cU{\mathcal{U}}
\newcommand\cV{\mathcal{V}}
\newcommand\cW{\mathcal{W}}

\newcommand\cY{\mathcal{Y}}

\newcommand\ve{\varepsilon}

\def\bA{\mathbf{A}}

\def\bC{\mathbf{C}}

\def\bE{\mathbf{E}}

\def\bH{\mathbf{H}}

\def\bP{\mathbf{P}}

\def\bR{\mathbf{R}}

\def\bZ{\mathbf{Z}}

\newcommand\del{\partial}

\newcommand{\GL}{\textup{GL}}
\newcommand{\Pic}{\textup{Pic}}

 \usepackage[all]{xy}
                \usepackage{amsthm,amsfonts,amsmath,amscd,amssymb,epsfig,verbatim,
}

\usepackage[parfill]{parskip}    
\usepackage{graphicx}
\usepackage{epstopdf}
\date{\today} 
 
\usepackage{graphicx}

\usepackage[dvipsnames,svgnames,x11names,hyperref]{xcolor}
\usepackage{url,graphicx,verbatim,amssymb,enumerate,stmaryrd}
\usepackage[colorlinks,citecolor=Bittersweet,linkcolor=Bittersweet,urlcolor=Bittersweet,filecolor=Bittersweet]{hyperref}
\usepackage{tikz} 
\usetikzlibrary{arrows,decorations.pathmorphing,backgrounds,fit,positioning,shapes.symbols,chains,calc}
\usepackage[all]{xy}
\xyoption{matrix}
\xyoption{arrow}
\usepackage[hmargin=3cm,vmargin=3cm]{geometry}
\usepackage{setspace,kantlipsum}
\tikzset{help lines/.style={step=#1cm,very thin, color=gray},
help lines/.default=.5} 
\usepackage{tikz-cd}
\usepackage{booktabs}
 
\theoremstyle{plain}
\newtheorem{theorem}{Theorem}[section]

\newtheorem{lemma}[theorem]{Lemma}
\newtheorem{proposition}[theorem]{Proposition}

\newtheorem{corollary}[theorem]{Corollary}

\newtheorem{cor}[theorem]{Corollary}

\theoremstyle{definition}
\newtheorem{definition}[theorem]{Definition}
\newtheorem{inductive step}[theorem]{Inductive step}

\newtheorem{inductive lemma}[theorem]{Inductive Lemma}
\theoremstyle{remark}
\newtheorem{example}[theorem]{Example}
\newtheorem{remark}[theorem]{Remark}
\newtheorem*{remark*}{Remark}
\newtheorem{notation}[theorem]{Notation}

\newtheorem{problem}[theorem]{Problem}
\newtheorem{question}[theorem]{Question}

\newtheorem*{example*}{Example}

\renewcommand{\mod}{\mathrm{mod}} 

\newcommand{\wh}{\widehat}

\newcommand{\inv}{^{-1}}

\newcommand{\bS}{{\bf S}}

\newcommand{\A}{{\bf A}}

\newcommand{\cL}{\mathcal{L}}

\def\Op{{\mathcal O}{\it p}\,}
\numberwithin{figure}{section}

\usepackage{stackengine}
\def\rlwd{.4pt}
\def\rlht{1.1pt}
\def\rhatvrule{\rule{\rlwd}{\rlht}}
\def\rhat#1{%
  \setbox0=\hbox{$#1$}%
  \stackon[0pt]{\stackon[1pt]{\ensuremath{#1}}{%
    \rhatvrule\kern\wd0\kern-\rlwd\kern-\rlwd\rhatvrule}}%
    {\rule{\wd0}{\rlwd}}%
}

\usepackage{tikz}
\usetikzlibrary{calc, arrows.meta}

\title{On arborealization, Maslov data, and lack thereof}

\author{Daniel  \'Alvarez-Gavela}
\address{Department of Mathematics \\ Brandeis University \\ Waltham, MA, 02453}
\email{dgavela@brandeis.edu}
\thanks{DA was partially supported by NSF grant DMS-1638352 and the Simons Foundation}

\author{Tim Large}
\email{tmjlarge@gmail.com}
\thanks{TL was supported by a Simons Junior Fellowship}

\author{Abigail Ward}
\address{Department of Pure Mathematics and Mathematical Sciences \\ University of Cambridge \\ Cambridge, UK, CB3 0WB}
\email{arw204@cam.ac.uk }
\thanks{AW was supported by NSF Grant DMS-2002183 and  UKRI Frontier Research Grant `Floer Theory Beyond Floer’ (grant number EP/X030660/1).}

\begin{document} 
\begin{abstract}

For a Weinstein manifold, we compare and contrast the properties of admitting an arboreal skeleton and admitting Maslov data. Both properties are implied by the existence of a polarization, for which a necessary condition is that the odd Chern classes be 2-torsion. In a similar spirit we establish cohomological obstructions to the existence of arboreal skeleta and to the existence of Maslov data, exhibiting concrete examples to illustrate their failure. For instance, we show that complements of smooth anti-canonical divisors in complex projective space may fail to admit either arboreal skeleta or Maslov data. We also exhibit an example of a Weinstein manifold which admits Maslov data but does not admit an arboreal skeleton. 

 \end{abstract}

\maketitle

 \onehalfspacing
 \tableofcontents
 
 \section{Introduction}
 
 \subsection{Main results}
 
 Let $X$ be a Weinstein manifold of dimension $2n$. Consider the following three properties: 
 
\begin{enumerate}
\item $X$ admits a polarization. \label{item1}
\item $X$ admits an arboreal skeleton. \label{item2}
\item $X$ admits Maslov data. \label{item3}
 \end{enumerate}
 
 That $X$ admits a polarization means that the tangent bundle contains a global field of Lagrangian planes, or equivalently that $TX \simeq E \otimes \bC$ as a complex vector bundle for some real vector bundle $E$, where $TX$ is endowed with the homotopically unique almost complex structure  compatible with the symplectic structure. Another equivalent condition is that the classifying map $X \to BU(n)$ for $TX$ lifts to a map $X \to BO(n)$. Given that $X$ has the topology of a CW complex of dimension at most $n$, this is equivalent to asking that the stabilization $X \to BU$ of the aforementioned classifying map lifts to a map $X \to BO$. 
 
 That $X$ admits Maslov data means that the Maslov obstruction, which is a certain classifying map 
$X \to B^2\Pic(\bS)$, is nullhomotopic. Here and throughout $\bS$ denotes the sphere spectrum. When the Maslov obstruction vanishes, one expects it to be possible to endow the moduli spaces of pseudo-holomorphic curves in $X$ with coherent orientations over $\bS$, see the discussion in Section \ref{sec: background} below. In \cite{NS} the vanishing of the Maslov obstruction appears as the necessary condition for defining a microlocal category for $X$ over $\bS$, invariant under Weinstein homotopy.  This is unsurprising because the microlocal category may be thought of as a sheaf-theoretic version of a Fukaya category of exact Lagrangians in $X$ with coefficients in $\bS$; see \cite{GPS} for a comparison over $\bZ$ in the polarizable case.
 
That $X$ admits an arboreal skeleton means that there exists a Weinstein structure in the Weinstein homotopy class of $X$ for which the skeleton is a stratified Lagrangian subset with singularities symplectomorphic to the standard models for arboreal singularities \cite{AGEN1}. Arboreal singularities were introduced by Nadler as a tractable class of Lagrangian singularities. They are indexed by finite rooted trees whose edges not adjacent to the root are labelled with $\pm$ signs. The microlocal category of an arboreal singularity is equivalent to modules over the quiver given by the corresponding directed tree \cite{Na1}. This gives a combinatorial description of the Fukaya or microlocal category in terms of the stratification of an arboreal skeleton.

We first consider the extent to which properties (\ref{item1}), (\ref{item2}) and  (\ref{item3}) fail in nature. 

Property (\ref{item1}), i.e., the existence of a polarization, is certainly a strong restriction on a Weinstein manifold $X$. Indeed, it is equivalent to asking that $TX$ is the complexification of a real vector bundle, and hence implies that 
\begin{equation}\label{eq:basic-obs}
2c_{2k+1}(X)=0
\end{equation}
 in $H^{4k+2}(X;\bZ)$ for all $k \geq 0$. Using this cohomological obstruction it is straightforward to exhibit examples of Weinstein manifolds which do not admit polarizations, for example the complement $X$ of a smooth divisor $D$ in the complex projective space $\bP^n$ only satisfies $2c_1(X)=0$ when the degree of $D$ divides $2(n+1)$.
 
Note that the cohomological condition $2c_1(X)=0$ is the familiar obstruction to the existence of $\bZ$-gradings in the Fukaya category, and is similarly known to obstruct the existence of Maslov data over $\bZ$, see Example \ref{ex:gradable} for further discussion. We will henceforth use the following terminology:

\begin{definition} We say that a Weinstein manifold $X$ is  {\em gradable}
 if  $2c_1(X)=0$ in $H^2(X;\bZ)$.\end{definition}

Property (\ref{item2}) also often fails to hold in nature. Indeed, we show that it implies the following cohomological obstruction:
\begin{theorem}\label{thm:obs-arb}
If a Weinstein manifold $X$ admits an arboreal skeleton, then \begin{equation}\label{eqn:obstruction} c_{2k+1}(X) = c_1(X)c_{2k}(X) \end{equation} in $H^{4k+2}(X;\bZ[\frac{1}{2}])$ for all $k \geq 1$.
\end{theorem}

In \cite{Na2} it is shown that it is always possible to locally deform any Legendrian singularity to an arboreal singularity while preserving the category of microlocal sheaves. So at least at the level of the microlocal category there are no obstructions to the local arborealization of a Lagrangian skeleton. However Theorem \ref{thm:obs-arb} asserts that there exist global obstructions. 

Property (\ref{item3}) is known to be obstructed by the gradability condition $2c_1(X)=0$ as noted above. We show that there exist further cohomological obstructions to property (\ref{item3}), namely:

\begin{theorem}\label{thm:obs-mas}
Let $p$ be an odd prime. If a Weinstein manifold $X$ admits Maslov data, and $c_i(X)=0$ in $H^*(X;\bZ/p)$ for $1 \leq i \leq p-1$, then also $c_p(X)=0$ in $H^*(X;\bZ/p)$. 
\end{theorem}

In fact for every odd prime $p$ we construct a characteristic class $m_p(X) \in H^{2p}(X;\bZ/p)$ which obstructs the vanishing of the Maslov obstruction and is a polynomial in the Chern classes $c_1,c_2,\ldots,c_p$ with nonzero $c_p$ coefficient. The consequence Theorem \ref{thm:obs-mas} will suffice for our sample applications. 

It is not hard to find examples of $X$ to which the above obstructions to the existence of arboreal skeleta and existence of Maslov data apply, and even examples which are gradable. For instance, let us consider examples of the form $X=\bP^n \setminus D$ for $D \subset \bP^n$ a smooth anti-canonical divisor, i.e., of degree $n+1$.

\begin{example}\label{ex:arb}
If $n>6$ is congruent to $2$ mod $6$, then the complement $X$ of a smooth anti-canonical divisor in $\bP^n$ does not admit an arboreal skeleton. For example the complement of a smooth anti-canonical divisor in $\bP^8$ does not admit an arboreal skeleton.
\end{example}

\begin{example}\label{ex:mas}
For any prime $p>2$ and any $k \geq 3$ not divisible by $p$, the complement $X$ of a smooth anti-canonical divisor in $\bP^{kp-1}$ does not admit Maslov data. For example the complement of a smooth anticanonical divisor in $\bP^{11}$ does not admit Maslov data. 
\end{example}

\begin{example} Taking $k=3$, $p=5$ we see from the previous two examples that the complement of a smooth anticanonical divisor in $\bP^{14}$ admits neither an arboreal skeleton nor Maslov data.
\end{example}

The spaces in Example \ref{ex:mas} admit wrapped Fukaya categories over $\bZ$, but the Maslov obstruction would seem to obstruct the refinement of these categories to categories over $\bS$.

\begin{question} Can the non-vanishing of the Maslov obstruction for the Weinstein manifolds of Example \ref{ex:mas}  be detected by their integer Fukaya categories?  \end{question}

\begin{question} What do the skeleta of the Weinstein manifolds of Example \ref{ex:arb} look like? In particular, what Lagrangian singularities beyond arboreal singularities are needed? \end{question}
 
We next consider the implications which hold between the properties (\ref{item1}), (\ref{item2}), and (\ref{item3}).

It was shown in work of the first author with Eliashberg and Nadler \cite{AGEN3} that existence of a polarization is equivalent to existence of a positive arboreal skeleton, hence in particular there holds the implication (\ref{item1}) $\Rightarrow$ (\ref{item2}). Note that this is consistent with the implication (\ref{eq:basic-obs})$ \Rightarrow $(\ref{eqn:obstruction}) for the cohomological obstructions to (\ref{item1}) and (\ref{item2}) considered above.

The implication (\ref{item1}) $\Rightarrow$ (\ref{item3}) is also known to hold. Indeed, recall that a choice of polarization furnishes a lift of the map $X \to BU$ to $X \to BO$ and therefore the composition $X \to BO \to BU \to B(U/O)$ is canonically nullhomotopic. The Maslov obstruction factors as a composition $X  \to BU \to B(U/O) \to B^2\Pic(\bS)$ and hence it is also canonically nullhomotopic. In the presence of a polarization, the problem of coherently orienting moduli spaces is much more tractable, and existent treatments of the construction of Floer homotopy types stay within this setting; see \cite{ADP, B, L, PS1, PS2}. 

Next, note that property (\ref{item2}) does not imply either of the other two properties. To see this we recall from work of Starkston \cite{St} that any 4-dimensional Weinstein manifold admits an arboreal skeleton, while the other two properties imply the gradability condition $2c_1(X)=0$ in $H^2(X;\bZ)$, which certainly does not hold for all 4-dimensional Weinstein manifolds $X$. 
 
We show that (\ref{item3}) $\not\Rightarrow$ (\ref{item2}) (and hence (\ref{item3}) $\not\Rightarrow$ (\ref{item1})) in Section \ref{sec:crazy-ex}, namely we derive the following consequence of Theorem \ref{thm:obs-arb};

\begin{corollary}\label{cor: ex Maslov} 
There exits a Weinstein manifold $X$ which admits Maslov data but which does not admit an arboreal skeleton.  \end{corollary}

Indeed, in Proposition \ref{prop: crazy-ex} we exhibit a Weinstein manifold homotopy equivalent to $S^6$ which admits Maslov data and such that $c_3(X) \neq 0$. It thus follows from Theorem \ref{thm:obs-arb} above that $X$ does not admit an arboreal skeleton, and certainly not a polarization. 

 The Weinstein manifold $X$ of Corollary \ref{cor: ex Maslov} can be thought of as a counterexample to a naive h-principle type expectation, namely that whenever the microlocal category of a Weinstein manifold can be defined algebraically (i.e., whenever the Maslov obstruction vanishes), one may in fact realize it geometrically in terms of the combinatorics of an arboreal skeleton.  
 
 \begin{question} Is there a larger class of tractable Lagrangian singularities for which such an h-principle might hold, or is any such expectation hopeless? \end{question}
 
In fact, the Weinstein manifold $X$ constructed in the proof of Corollary \ref{cor: ex Maslov} is obtained via the existence h-principle for Weinstein structures and hence is a priori flexible. Thus in this case any obstruction to arborealizability is invisible from the point of view of the Fukaya category of $X$, which is necessarily trivial if $X$ is flexible. In Remark \ref{rem: final} we explain how one may modify the construction to obtain a non-flexible $X$ as in Corollary \ref{cor: ex Maslov}.

\begin{remark}
It is currently not known whether an h-principle holds for the existence of (not necessarily positive) arboreal skeleta for Weinstein manifolds. Said differently, it is not known whether property (\ref{item2}) can be characterized in homotopy theoretic terms, similarly to how existence of a positive arboreal skeleton is equivalent to existence of a polarization. The notion of an $\cI$-structure introduced in Section \ref{sec: I-str} is the best known candidate, but it is not known whether this necessary condition is also sufficient. New ideas are needed to go beyond the polarizable case treated in \cite{AGEN3}.
\end{remark}

\subsection{The Maslov obstruction}\label{sec: background} Let us for a moment give an impressionistic account of why the (non)existence of polarizations or Maslov data arises in symplectic geometry. Suppose $X$ is an almost complex manifold with a distinguished base-point, and that we have a family $M$ of based smooth maps  $\bP^1 \to X$. 
This family has an evaluation map ${ev}: M \times S^2 \to X$ and a complex vector bundle ${ev}^*TX$ from pulling back the tangent bundle of X. Consider the $\overline{\del}$ operator for sections of this bundle on the $S^2$ fibers of the projection $M \times S^2 \to M$: its index bundle defines a virtual complex vector bundle $\mathrm{ind}(\overline{\del})$ on $M$. In the event that $X$ is symplectic and $M$ is in fact a smooth moduli space of one-pointed pseudoholomorphic maps, this index bundle coincides with the class of the tangent bundle $[TM]$. Regardless of this, the index bundle has a very simple characterization following Atiyah’s proof of Bott periodicity \cite{At}. The map $M \times S^2 \to BU$ classifying ${ev}^* TX$ yields a map
\[M \to \Omega^2 BU\]
and composing this with the Bott equivalence $\Omega^2 BU \simeq \bZ \times BU$ gives a map $M \to \bZ \times BU$ which coincides up to homotopy with the classifying map of the index bundle  $\mathrm{ind}(\overline{\del})$. 

Floer- and Gromov-Witten-theoretic invariants valued in a ring $R$ could then be defined by ``counting" $M$ in some sense. This requires that $TM$ is oriented over $R$: this is the classical notion for ordinary rings, but for ring spectra $R$, this means that the composite map
\[ M \to  \bZ \times BU \to \Pic(R)\]
should be nullhomotopic (see \cite{An}).

Moreover, the interest in Floer and Gromov-Witten theory comes through the algebraic structures these invariants form; this requires that, when variations of the above construction are done with more base-points, constraints, or Lagrangian boundary conditions, these orientations must be coherent under the appropriate gluing operation.

The primary way in which coherent orientations can be constructed in practice is to exhibit a structure on the tangent bundle $TX$ of $X$ itself, which induces a certain structure on the tangent bundle $TM$ of a transversely cut out moduli space $M$ through the identification of $[TM]$ with the index bundle. For example, if the map $X \to BU$ was null-homotopic, so would be the map $M \times S^2 \to BU$ and thus the map $M \to \Omega^2 BU \simeq \bZ \times BU$ which classifies the index bundle. However in practice the obstruction is looser: to construct Floer theoretic invariants of $X$ valued in any ring spectrum, it is enough for the tangent bundles of $M$ to be (coherently) trivial as stable sphere bundles, a much weaker condition than being trivial as a stable complex vector bundle. This amounts to saying that the composite map
\[ X \to BU \to B^2(\bZ \times BU) \to B^2 \Pic(\bS)\]
is null-homotopic, where the map $BU \to B^2(\bZ \times BU)$ is complex Bott periodicity, $\bS$ is the sphere spectrum, and the map $B^2(\bZ \times BU) \to B^2 \Pic(\bS)$ is the three-fold delooping of the complex J-homomorphism. However it is often more convenient to recast this using real Bott periodicity; this map is homotopic to the composite
\[ X \to BU \to B(U/O) \to B^2(\bZ \times BO) \to B^2 \Pic(\bS) \]
where the maps $B(U/O) \to B^2(\bZ \times BO)$ and $B^2(\bZ \times BO) \to B^2 \Pic(\bS)$ are real Bott periodicity and the three-fold delooping of the real J-homomorphism respectively. We call the map $X \to B^2\Pic(\bS)$ the (pseudo-holomorphic) {\em Maslov obstruction}.

The role of Picard groups of ring spectra in this story was first noted in \cite{Lu}, but for ordinary rings this formulation of index theory for $\overline{\del}$ operators has always underpinned the ability to define Floer homology with gradings and signs (see, e.g., \cite{FlHo} and \cite[Section II.11]{Seidel}).

In practice, to define sphere-spectrum valued Fukaya categories, it’s convenient to work with the stronger assumption that $X$ admits a polarization, which implies $X \to B(U/O)$ is null-homotopic and thus so is the full composition $X \to B^2 \Pic(\bS)$, since the additional geometric structure simplifies the task of defining coherent gluing maps.

From the microlocal viewpoint, given a smooth closed exact Lagrangian submanifold $L$ of a cotangent bundle $T^*M$, the microlocal sheaf theory of Kashiwara–Schapira naturally assigns a local system of categories on $L$ with fiber equivalent to the category of spectra \cite{JiTr}. This is the spectral upgrading to the microlocal approach to the Fukaya category of cotangent bundles \cite{NZ}, \cite{G}. The foundational work of Jin \cite{J1,J2} shows that the classifying map $L \to B\Pic(\bS)$ of this local system of categories factors as a composition $L \to U/O \to B \Pic(\bS)$ where $L \to U/O$ is the Gauss map and $U/O \to B\Pic(\bS)$ classifies a local system of categories on $U/O$. Further, Jin identifies the map $U/O \to B\Pic(\bS)$ as the composition of the real Bott periodicity map $U/O \to B(\bZ \times BO)$ and the two-fold delooping of the real $J$-homomorphism $B(\bZ \times BO) \to B\Pic(\bS)$.

One loop down, the (microlocal) Maslov obstruction of a Weinstein manifold $X$ is a map $X \to B^2 \Pic(\bS)$ which factors as a composition $X \to BU \to B(U/O) \to B^2\Pic(\bS)$, where $X \to BU$ is the classifying map of $TX$ and $B(U/O) \to B^2\Pic(\bS)$ is a universal map which classifies a local system of 2-categories on $B(U/O)$. This universal map is a delooping of the microlocal classifying map $U/O \to B\Pic(\bS)$. When the Maslov obstruction of a Weinstein manifold $X$ vanishes, i.e. is nullhomotopic, Nadler and Shende show \cite{NS}  that one may define a microlocal category for $X$ over $\bS$, invariant under Weinstein homotopy.

It follows from the work of Jin cited above that the microlocal Maslov obstruction of a Weintein manifold $X$ also factors as
$$ X \to BU \to B(U/O)  \to  B^2(\bZ \times BO) \to B^2\Pic(\bS).$$
where the maps $B(U/O) \to B^2(\bZ \times BO)$ and $B^2(\bZ \times BO) \to B^2 \Pic(\bS)$ are real Bott periodicity and the three-fold delooping of the real J-homomorphism respectively.

To verify that the pseudo-holomorphic and the microlocal Maslov obstructions coincide, one would need to compare the model for complex Bott periodicity $\Omega^2BU \simeq \bZ \times BU$ from Atiyah \cite{At} with the model for real Bott periodicity $\Omega (U/O)  \simeq \bZ \times BO$ used by Jin \cite{J1,J2}. It would be sufficient to know that the resulting diagram
\[  \xymatrix{  \Omega^2BU \ar[d] \ar[r] & \bZ \times BU   \ar[d] \\ \ \Omega^2B(U/O) \ar[r]&  \bZ \times BO } \] 
is homotopy commutative at the level of $\bE_2$-maps. One might expect this to be the case, even at the level of infinite loop maps, but we do not know of such a comparison in the literature.

\begin{remark}  In any case, given any realization of Atiyah's Bott equivalence $\Omega^2BU \simeq \bZ \times BU$ as an $\bE_2$-map such that the resulting map $BU \to B^2\Pic(\bS)$ factors up to homotopy as the composition of the map $BU \to B(U/O)$, some model  for real Bott periodicity $B(U/O) \simeq B^2(\bZ \times BO)$, and the three-fold delooping of the real J-homomorphism $B^2(\bZ \times BO) \to B^2\Pic(\bS)$, this map will only differ from the microlocal Maslov obstruction by the insertion of some homotopy automorphism $B(U/O) \simeq B(U/O)$. The cohomological obstruction to the vanishing of the Maslov obstruction that we construct in Theorem \ref{thm: char class} is a class in $H^{2p}(BU;\bZ/p)$ which is the pullback of a class in $H^{2p}(B(U/O);\bZ/p)$. If we replace this class with its pullback under any homotopy automorphism $B(U/O) \simeq B(U/O)$ we obtain a class in $H^{2p}(BU;\bZ/p)$ which also satisfies the properties of Theorem \ref{thm: char class}. Hence the conclusion of Theorem \ref{thm: char class} holds for either Maslov obstruction $X \to B^2\Pic(\bS)$ independently of whether the homotopy automorphism $B(U/O) \simeq B(U/O)$ which relates them is the identity, and so does Theorem \ref{thm:obs-mas}.  \end{remark}

\subsection{Structure of the article}

In Section \ref{sec: I-str} we introduce the notion of an $\cI$-structure on a complex vector bundle and establish a cohomological obstruction to existence of an $\cI$-structure. In Section \ref{sec: arb} we show that existence of an arboreal skeleton implies existence of an $\cI$-structure. In Section \ref{sec: Maslov} we establish a cohomological obstruction to existence of Maslov data and construct a Weinstein manifold which admits Maslov data but is not arborealizable.

\subsection{Acknowledgements}
We are grateful to John Pardon, who suggested the notion of an $\cI$-structure as a nontrivial obstruction to the existence of an arboreal skeleton, and to S$\o$ren Galatius who indicated how one might obtain cohomological obstructions to the existence of $\cI$-structures. We are also grateful to Yasha Eliashberg and David Nadler for their perspectives on Weinstein manifolds, whose influence on this article is apparent. We are also grateful to Xin Jin and Oleg Lazarev for useful conversations on various topics related to this work. We are also grateful to Paul Seidel for his support at MIT, where many conversations related to the ideas in this article took place, and to Ivan Smith for many helpful conversations, as well as comments on a draft. We are also grateful to Kenneth Blakey and Ceyhun Elmacioglu for pointing out mistakes in an earlier version of the article.

\section{$\cI$-structures} \label{sec: I-str}
\subsection{Definitions}

Let $(X,\omega)$ be a symplectic manifold of dimension $2n$. Recall that $X$ is said to be polarizable if there exists a global Lagrangian distribution in $TX$, i.e., if the  Lagrangian Grassmanian bundle of $TX$ admits a global section. Equivalently, we can state the condition in terms of the classifying map $X \to BU(n)$ for $TX$ equipped with the Hermitian structure associated to the contractible choice of an $\omega$-compatible almost complex structure. Namely, $X$ is polarizable if and only if there exists a lift 
\[  \xymatrix{  & BO(n) \ar[d]  \\ X \ar@{.>}[ur] \ar[r]^{[TX]} & BU(n) }. \]
We can relax this condition by defining the notion of an $\cI$-structure on a Hermitian vector bundle. This notion was suggested by J. Pardon to the authors of \cite{AGEN1} as a nontrivial obstruction to the existence of an arboreal skeleton.

\begin{definition}
Let $\cI(n)$ be given by the homotopy pushout diagram 
\[ \xymatrix{  BO(n-1) \times BO(1) \ar[d]^{f_0} \ar[r]^{f_1}& BO(n-1) \times BU(1)  \ar[d] \\ BO(n) \ar[r] & \cI(n) } \]
where the map $f_0: BO(n-1) \times BO(1) \to BO(n)$ is given by the direct sum and the map $f_1: BO (n-1) \times BO(1) \to BO(n-1) \times BU(1)$ is given by complexification on the second factor. 
\end{definition}

\begin{definition} We denote by $g:\cI(n) \to BU(n)$ the natural map inherited from the maps $BO(n) \to BU(n)$ and $BO(n-1) \times BU(1) \to BU(n)$. \end{definition}

\begin{definition} 
Let $X$ be a paracompact topological space and $\cW \to X$ a Hermitian vector bundle of complex rank $n$ with classifying map $[\cW]:X \to BU(n)$. An {\em $\cI$-structure} on $\cW$ is defined to be a lift 
\[  \xymatrix{  & \cI(n) \ar[d]^{g}  \\ X \ar@{.>}[ur] \ar[r]^{[\cW]} & BU(n)} \] up to homotopy.
\end{definition}

\begin{figure}[h]
\includegraphics[scale=0.15]{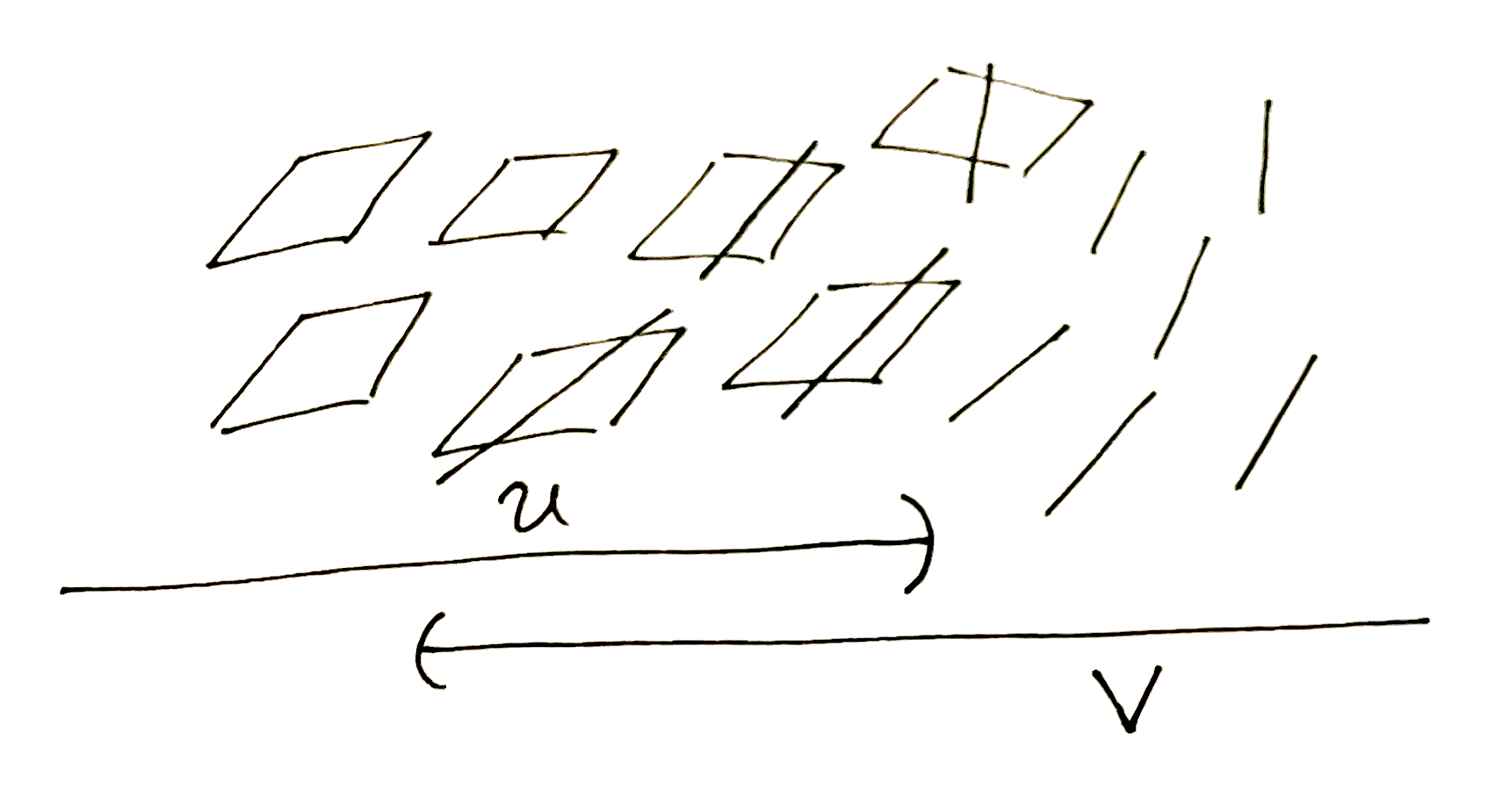}
\caption{An illustration of an $\cI$-structure subordinate to an open cover $U \cup V$.}
\label{fig:I-str}
\end{figure}

It will be convenient to have a more concrete description of this notion in terms of an open cover for $X$, see Figure \ref{fig:I-str}.

\begin{definition}
Let $X$ be a paracompact topological space and $\cW \to X$ be a Hermitian vector bundle of complex rank $n$. An $\cI$-structure {\em subordinate} to an open cover $ X = U \cup V$ consists of the following data:
\begin{enumerate} 
\item Over $U$, a dimension $n$ Lagrangian distribution $\cF \subset \cW|_U$;
\item Over $V$  a dimension $(n-1)$ isotropic distribution $\cG \subset \cW|_{V}$;
\item Compatibility over $U \cap V$ in the sense that 
\[ \cG_p \subset \cF_p \text { for all } p \in U \cap V.\]
\end{enumerate} 
\end{definition}

\begin{lemma}\label{lem: open cover} Let $\cW \to X$ be a Hermitian vector bundle. Then $\cW$ admits a $\cI$-structure if and only if there exists an open cover $ X = U \cup V$ with a subordinate $\cI$-structure.
\end{lemma}

\begin{proof}  Take the following realization of the homotopy colimit: 
\[ \cI(n) = \frac{BO(n) \sqcup( BO(n-1) \times BO (1) \times [0,1] )\sqcup BO(n -1) \times BU(1)}{ ((x,0) \sim f_0(x), (x,1) \sim f_1(x))}. \]
This space has a cover satisfying the conditions of the theorem: letting $P: \cI(n) \to [0,1]$ be the projection, we take
\[U_{\cI} = P^{-1}( [0,1)), \,\,\,\, V_{\cI} = P^{-1}((0,1]).\] There are natural maps $\tilde{f}_0: U_{\cI} \to BO(n)$ and $\tilde{f}_1: V_{\cI} \to BO(n-1) \times BU(1)$, and we take our distributions to be 
\[ \cF_{\cI} = (\tilde{f}_0)^* E O(n), \,\,\,\, \cG_{\cI} = (\tilde{f}_1)^* E O(n-1), \]
where both $E O(n)$ and $E O(n-1)$ are viewed as sub-bundles of $E U(n)$ under the complexification maps $g$  to $BU(n)$. Note $\cG_{\cI} \subset \cF_{\cI}$ on $U_{\cI} \cap V_{\cI}$ by design.

Now assume that $\cW \to X$ has an $\cI$-structure classified by a map $p: X \to \cI(n)$, so we can write $\cW \cong p^*g^*EU(n).$ Then taking $U= p\inv(U_{\cI}), V=p\inv(V_{\cI})$, $\cF = p^*\cF_{\cI}$ and $\cG=p^*(\cG_{\cI})$ yields our desired cover. 

Conversely, assume that $X$ has such a cover. Then the space $X$ is homotopy equivalent to the homotopy colimit of the inclusions $i_U$ and $i_V$ of $U \cap V$ into $U$ and $V$ respectively, so we can consider instead the thickened space 
\[Y = \frac{U \cup ((U \cup V) \times [-1,2]) \cup V}{((x,-1) \sim i_U(x), (x,2) \sim i_V(x))},  \]
with the vector bundle $\pi^*\cW$, where $\pi: Y \xrightarrow{\simeq} X$ is the natural map. We  now explicitly construct a map $Y \to \cI(n)$ lifting the map classifying $\pi^*\cW$.

Taking orthogonal complements, we can choose a complex line bundle $\cL \subset \cW$ on V such that \[ (\cG \otimes \bC ) \oplus \cL = \cW|_V; \] then after taking $\cL_{\bR} = \cL|_{U \cap V} \cap \cF|_{U \cap V}$, we can write \[ \cF|_{U \cap V} = \cG|_{U \cap V} \oplus \cL_{\bR} .\]  Let $p_0: U \to BO(n)$ be a map classifying $\cF$, let $p_1: V \to BO(n - 1) \times BU(1)$ classify the pair $(\cG, \cL)$ and let $p_2: U \cap V \to BO(n-1) \times BO(1)$ classify $\cG \oplus \cL_{\bR}$. 
We then define a map $p: Y \to \cI(n)$ by gluing the maps
\begin{align*} &p_0: U  \to BO(n);\\
		    &h_0: U \cap V \times [-1,0] \to BO(n); \\
		    &(p_2, \mathrm{id}): U \cap V \times [0,1] \to BO(n-1) \times BO(1) \times [0,1] ; \\
		    &h_1  : U \cap V \times [1,2] \to BO(n-1) \times BU(1),\\
		   & p_1: V \to BO(n-1) \times BU(1) \end{align*}
where $h_0$ is a homotopy between $p_0|_{U \cap V}$ and $f_0 \circ p_2$ (which necessarily exists since both maps classify $\cF$), and $h_1$ is a homotopy between $p_1|_{U \cap V}$ and $f_1 \circ p_2$ (which similarly exists since both maps classify the pair $(\cG, \cL)$). By construction, the complexification of every map in sight classifies the restriction of the Hermitian vector bundle $\pi^*\cW$ and hence $(g \circ p)^* E U(n)$ is isomorphic to $\cW$. We conclude that $\cW$ has an $\cI$-structure. \end{proof}
		    		    
\begin{remark} There are natural stabilization maps $\cI(n) \to \cI(n+1)$ induced by the stabilization maps $BO(n) \to BO(n+1)$, $E \mapsto \bR \oplus E$. Given a Hermitian vector bundle $\cW \to X$ with an $\cI$-structure $(U, V, \cF, \cG)$ classified by a map $X \to \cI(n)$,  
the composition
\[ X \to \cI(n) \to \cI(n+1) \]
classifies the $\cI$-structure  $(U, V,  \underline{\bR} \oplus \cF ,  \underline{\bR} \oplus \cG)$ on the Hermitian vector bundle $\underline{\bC}  \oplus \cW \to X$. 
We define $\cI = \lim_{n} \cI(n)$. There is a map $\cI \to BU$ and a homotopy pushout diagram 
\[ \xymatrix{  BO \times BO(1) \ar[d] \ar[r] & BO \times BU(1)  \ar[d] \\ BO \ar[r] & \cI }. \] 
We may talk of {\em stable} $\cI$-structures on $\cW$ as (homotopy classes of) lifts:
\[  \xymatrix{  & \cI \ar[d]^{g}  \\ X \ar@{.>}[ur] \ar[r]^{[\cW]} & BU } \] 
\end{remark}

\begin{remark}\label{rem: lift} Since the inclusion $U(n) \subset \GL_n(\bC)$ is a homotopy equivalence, every complex vector bundle has a homotopically unique Hermitian structure; thus we may unambiguously speak of an $\cI$-structure on a complex vector bundle without having to specify the Hermitian structure. Similarly, since $U(n) \subset \text{Sp}_{2n}(\bR)$ is a homotopy equivalence, we may also speak of an $\cI$-structure on a real symplectic vector bundle without having to specify the Hermitian structure. 
\end{remark}

\subsection{A cohomological obstruction}

For a complex vector bundle, existence of an $\cI$-structure is not as strong a condition as being the complexification of a real vector bundle, but it is nonetheless quite restrictive. Concretely, we will show that existence of an $\cI$-structure implies the following restriction on the Chern classes of a complex vector bundle.

\begin{proposition}\label{prop:lift-obs-main}Let $\cW \to X$ be a complex vector bundle. 
 If $\cW$ admits an $\cI$-structure, then for all $k \geq 1$,  \begin{equation} \label{eqn:Istruc} c_1(\cW) c_{2k}(\cW)- c_{2k+1}(\cW)=0 \text { in  } H^{4k+2}(X; \bZ[\tfrac{1}{2}]). \end{equation} 
\end{proposition}
\begin{proof}  It suffices to obstruct the stabilized problem of lifting the classifying map $X \to BU$ to a map $X \to \cI$. Since $\cI$ is a homotopy pushout, there is a Mayer-Vietoris sequence  
\begin{multline}   \cdots  \to H^*(\cI; R) \to H^{*} (BO; R) \oplus H^*(BO \times BU(1); R) \to H^*(BO \times BO(1); R)  \to \cdots \end{multline}
for any coefficient ring $R$. Now take $R$ to be $\bZ[\frac{1}{2}]$. We have 
 \begin{align*} H^*(BU;R) &\simeq R[c_1,c_2,\ldots,]; \\ 
 	       H^*(BO; R) &\simeq R[p_1, p_2, \ldots];  \\
	       H^*(BO \times BU(1);R) &\simeq R [ p_1, p_2, \ldots ] \otimes R[c_1]; 
	     \end{align*}
	    where $|c_i|=2i$ and  $|p_i|=4i$. (Strictly speaking, the second $c_1$ appearing is the pullback of $c_1 \in H^*(BU; R)$ under the map $BU(1) \to BU$, but this notational ambiguity does not affect the following calculation.)  Moreover, the map $BO \times BO(1) \to BO$ induces an isomorphism on cohomology with coefficients in $R$, so using the long exact sequence we can identify
 \[H^*(\cI; R)\cong \ker \left( H^{*} (BO; R) \oplus H^*(BO \times BU(1); R) \to H^*(BO; R)\right) \cong H^*(BO \times BU(1); R), \]
and under this isomomorphism, the map $H^*(BU; R) \to H^*(\cI; R)$ is identified with the map $H^*(BU; R) \to H^*(BO \times BU(1); R)$. We now analyze this map.

Suppose that $\cV \to X$ is a real vector bundle and $\cL \to X$ a complex line bundle. Then the total Chern class of $(\cV \otimes \bC) \oplus \cL$ is given by
\[ c\left( (\cV \otimes \bC) \oplus \cL \right)   = c(\cV \otimes \bC) \cdot c(\cL)  = (1-p_1(\cV)+p_2(\cV) - \cdots ) \cdot (1+c_1(\cL)) \]
\[= (1 \cdot1) + 1 \cdot c_1(\cL) - p_1(\cV) \cdot 1 -   p_1(\cV) \cdot c_1(\cL)  + p_2(\cV) \cdot 1  + p_2(\cV) \cdot c_1(\cL) + \cdots
\]

So we can compute the effect of the map $H^*(BU;R) \to H^*(BO \times BU(1) ;R) $ in terms of the basis given by Chern classes:
\[ 1 \mapsto 1, \quad  c_1  \mapsto 1 \cdot c_1 , \quad  c_2 \mapsto - p_1 \cdot 1 \] \[ c_3 \mapsto - p_1 \cdot c_1,  \quad c_2 c_1  \mapsto -p_1 \cdot c_1, \quad c_1^3 \mapsto 1 \cdot c_1^3, \qquad \cdots  \]
Note in particular that $c_{2k+1} \mapsto (-1)^{k} p_{k} \cdot c_1 $ and $c_{2k} \mapsto (-1)^{k} p_k \cdot 1$ for all $k \geq 1$. Hence we have that $c_{2k}c_1 -c_{2k+1}$ is in the kernel of the map $H^*(BU;R) \to H^*(BO \times BU(1) ;R)$. Since the existence of an $\cI$-structure on $\cW$ implies that the induced  map $[\cW]^*: H^*(BU; R) \to H^*(X; R)$ factors through the map to $H^*(BO \times BU(1) ;R)$, we conclude that Equation \eqref{eqn:Istruc}  holds.  \end{proof}

\subsection{Dual $\cI$-structures}
It will also be useful to have the following equivalent formulation of the notion of $\cI$-structure:

\begin{definition}
Let $\cW \to X$ be a symplectic vector bundle. A {\em dual $\cI$-structure} subordinate to an open cover $ X = U \cup V$ consists of the following data:
\begin{enumerate} 
\item Over $U$, a dimension $n$ Lagrangian distribution $\cF \subset \cW|_U$;
\item Over $V$  a dimension $(n+1)$ co-isotropic distribution $\cH \subset \cW|_{V}$;
\item Compatibility over $U \cap V$ in the sense that 
\[ \cH_p \pitchfork \cF_p \text { for all } p \in U \cap V.\]
\end{enumerate}  
 \end{definition}

In the following lemma the superscript indicates (real) dimension.

\begin{lemma}\label{lem:trans-contr}
Given a symplectic vector space $(W^{2n},\omega)$, and given a co-isotropic subspace $H^{n+k} \subset W$, the space $\Lambda_H(W)$ of isotropic subspaces $G^{n-k} \subset W$ which are transverse to $H$ is contractible. Similarly, given an isotropic subspace $G^k \subset W$, the space $\Lambda_G(W)$ of co-isotropic subspaces $H^{n+k} \subset W$ which are transverse to $G$ is contractible. 
\end{lemma}

\begin{proof}
The critical case $k=0$ in which $H$ is Lagrangian is standard: we identify $\Lambda_H(W)$ with the space of quadratic forms on any Lagrangian subspace transverse to $H$ and deduce that $\Lambda_H(W)$ is contractible.

Now, for any isotropic subspace $V$ of dimension $k$,  define a subset $\Omega_{W}(V)$ of the real Grassmanian of $2k$-dimensional subspaces of $W$ by
\[ \Omega_{V}(W) =\{K \in \mathrm{Gr}_{2k}(W) \mid V \subset K \text{ and } K \text{ is symplectic}\}. \] 

For a co-isotropic subspace $H \subset W$, we claim that there is a homotopy equivalence 
\begin{equation} \label{eqn:GH} \Lambda_H(W) \xrightarrow{\simeq} \Omega_{H^{\omega_{\perp}}}(W), \end{equation}
given by the map $ G \mapsto G + H^{\omega_\perp}$, and that the latter space is contractible. 

We first verify that this map is well-defined: indeed, note that given any $G \in \Lambda_{H}(W)$, we have  $G \cap H^{\omega_{\perp}} \subset G \cap H = 0$, and 
\[ (G + H^{\omega_{\perp}}) \cap (G + H^{\omega_{\perp}})^{\omega_\perp} = (G + H^{\omega_{\perp}}) \cap H \cap G^{\omega_{\perp}}= H^{\omega_{\perp}} \cap G^{\omega_{\perp}}=(H+G)^{\omega_{\perp}}  = 0. \]

The map of Equation \eqref{eqn:GH} is a Serre fibration. To see this, let $\text{Sp}(W,H,\omega)$ denote the group of linear symplectic automorphisms of $(W,\omega)$ which preserve $H$. For any fixed $K_0 \in \Omega_{H^{\omega_{\perp}}}(W)$ the map $\varphi \mapsto \varphi(K_0)$ defines a Serre fibration $\text{Sp}(W,H,\omega) \to \Omega_{H^{\omega_{\perp}}}(W)$. For any $G_0 \in \Lambda_H(W)$ such that $G_0+H^{\omega_\perp} = K_0$ we may factor this Serre fibration as the composition of the map $\text{Sp}(W,\omega) \to \Lambda_H(W)$, $\varphi \mapsto \varphi(G_0)$ and the map of Equation  \eqref{eqn:GH}, hence this latter map is also a Serre fibration.

The Serre fibration of Equation \eqref{eqn:GH} has fiber $\Lambda_{H^{\omega_{\perp}}}(K)$ over a point $K$. Hence by the $k=0$ case the fibers are contractible and thus the map is a homotopy equivalence as claimed.

It remains to show that $\Omega_{V}(W)$ is contractible for any co-iostropic subspace $V$. Suppose that we have a family $K^z \in \Omega_V(W)$ parametrized by $z \in \partial D^r$ for $D^r$ the unit disk of dimension $r$. We can choose compatible complex structures $J_1^z$ and $J_2^z$ on the two symplectic vector bundles $(K^z, \omega|_{K^z}) \to \partial D^r$ and $((K^z)^{\omega_{\perp}}, \omega|_{(K^z)^{\omega_{\perp}}}) \to \partial D^r$, respectively. Then $J_1^z \oplus J_2^z$ defines an $\omega$-compatible complex structure on the trivial vector bundle $\underline{W} \to \partial D^r$. The space of all $\omega$-compatible complex structures on $W$ is contractible and thus $J_1^z \oplus J_2^z$ extends to a complex structure $J^z$ on $\underline{W} \to D^r$. This allows us to extend $K^z$ over $D^r$ by taking 
\[K^z=V+J^zV \in \Omega_V(W), \] so we must have $\pi_r (\Omega_V (W))= 0$. Since this holds for any $r \geq 0$ we conclude that $\Omega_{V}(W)$ is contractible. 

The above discussion proves the statement for co-isotropic subspaces. To prove the corresponding statement for isotropic subspaces, we appeal to the fact that for any co-isotropic subspace $H$, the map $G \mapsto G^{\omega_{\perp}}$ realizes a diffeomorphism $\Lambda_H(W) \simeq   \Lambda_{H^\omega_{\perp}} (W).$ \end{proof}

\begin{proposition}\label{prop: dual}
Let $\cW \to X$ be a symplectic vector bundle. The data of an $\cI$-structure and of a dual $\cI$-structure subordinate to a given open cover $ X = U \cup V$ with fixed Lagrangian distribution $\cF$ on $U$ are homotopically equivalent. 
\end{proposition}

\begin{proof}
Given an $\cI$-structure $(\cF,\cG)$ subordinate to $(U,V)$, by Lemma \ref{lem:trans-contr} we may find an $(n+1)$-dimensional co-isotropic distribution $\cH$ on $V$ which is transverse to $\cG$, and moreover $\cH$ is unique up to contractible choice. Along $U \cap V$ we have $\cG \subset \cF$ and $\cG \cap \cH = 0$, hence $\dim \cH \cap \cF =1$ and consequently $\cH \pitchfork \cF$ as required. 

Conversely, given a dual $\cI$-structure $(\cF,\cH)$ subordinate to $(U,V)$, we have $\cL= \cH \cap \cF$ a rank 1 sub-bundle of $\cF$ along $U \cap V$ and hence there is a contractible choice of a hyperplane field $\cG \subset \cF$ on $U \cap V$ transverse to $\cL$. This $\cG$ is in particular an isotropic distribution on $U \cap V$ transverse to $\cH$ and hence by Lemma \ref{lem:trans-contr} can be extended to an isotropic distribution on $V$ transverse to $\cH$, uniquely up to contractible choice. \end{proof}

\section{Arboreal skeleta}\label{sec: arb}

\subsection{Weinstein manifolds}

\begin{definition}
A Weinstein structure on a manifold $X$ consists of a pair $(\omega, V)$ such that
\begin{enumerate}
\item $\omega$ is a symplectic form.
\item $V$ is a Liouville vector field, i.e., $\mathcal{L}_V \omega = \omega$.
\item there exists $\phi:V \to \bR$ a smooth exhausting function (proper and bounded below) such that the pair $(V,\phi)$ is gradient like in the sense of Cieliebak \cite{Cie24}. 
\end{enumerate}
\end{definition}

Sometimes it is simpler to restrict one's attention to Weinstein manifolds of finite type, though this won't be necessary for our discussion. A Weinstein manifold $X$ is said to be of {\em finite type} if there exists a compact manifold with boundary $W \subset V$, called a defining Weinstein domain, such that $V$ is outwards pointing along $\partial W$ and $\phi$ has no critical points outside of $W$. If a Weinstein manifold is of finite type then it is the completion of Weinstein domain in a standard way.

\begin{definition}
The {\em skeleton} of a Weinstein manifold is the attractor $\cK= \bigcap_{t>0} \varphi^{-t}(X)$ of the negative flow $\varphi^t$ of the Liouville vector field $V$.
\end{definition}

\begin{notation} For $A$ a subset of a topological space, we use the Gromov notation of letting $\Op(A)$ denote an arbitrarily small but unspecified open subset containing $A$. \end{notation}

\begin{lemma}\label{lem:contract}
Let $X$ be a Weinstein manifold and $\cK \subset X$ its skeleton. Given any open subset $\Omega \supset \cK$, there exists a smooth function $f:X \to [0,\infty)$ such that $\varphi^{-f(x)}(x) \in \Omega$ for all $x \in X$. 
\end{lemma}

\begin{proof}
It suffices to show that given $x \in X$ there exists $T>0$ such that for all $t>T$ we have $\varphi^{-t}(x) \in \Omega$; the construction of a suitable $f$ is then straightforward using a suitable partition of unity. Suppose then for contradiction that $x \in X$ is such that there exists a sequence $t_1<t_2<\cdots <t_k \to \infty$ such that $\varphi^{-t_i}(x) \not \in \Omega$. Since the points $\varphi^{-t_i}(x)$ stay in the compact subset $\phi^{-1}( -\infty, \phi(x) ]$, by passing to a subsequence we may assume that they converge to some limit $y \in X$. On the one hand we must have $y \not \in \cK$, since this would imply that $\varphi^{-t_i}(x) \in \Omega$ for $i$ large enough. On the other hand we must have $V(y)=0$, since otherwise the negative flow of $V$ would take points near $y$ outside of a neighborhood of $y$, never to come back. But $V(y)=0$ implies that $y \in \cK$, which is a contradiction. \end{proof}

\begin{remark}
If $X$ is of finite type and has a defining Weinstein domain $W \subset X$, then for any open subset $\Omega \supset \cK$ there exists a uniform $T_W>0$ such that $\varphi^{-T_W}(W) \subset \Omega$. 
\end{remark}

\begin{lemma} Existence of an $\cI$-structure on a Weinstein manifold $X$ is equivalent to existence of an $\cI$-structure   on $TX|_\cK$, the restriction of $TX$ to the skeleton $\cK \subset X$. 
\end{lemma}
\begin{proof}  Certainly any $\cI$-structure on $X$ restricts to an $\cI$-structure on $TX|_\cK$. Conversely, suppose that there exists a $\cI$-structure on $TX|_\cK$. Then there exists a dual $\cI$-structure $(\cF,\cH)$ subordinate to some cover $\cK=U \cup V$ by relatively open subsets. Choose any extension of $\cF$ (resp. $\cH$) to a Lagrangian (resp. co-isotropic) distribution on a neighborhood of $U$ (resp. a neighborhood of $V$) in $X$. Since  the condition $\cF \pitchfork \cH$ is open, for sufficiently small open neighborhoods $U'\subset X$ of $U$ and $V' \subset X$ of $V$, the pair $(\cF,\cH)$ will restrict to a dual $\cI$-structure on the open subset $\Omega=U' \cup V' \subset X$. Therefore there exists an $\cI$-structure on $\Omega \supset \cK$, i.e., the map $\Omega \to BU(n)$ classifying $T\Omega$ lifts to a map $\Omega \to \cI(n)$. 

By Lemma \ref{lem:contract}, there exists a function $f:X \to [0,\infty)$ such that $\varphi^{-f(x)}(x) \in \Omega$ for all $x \in X$. Let $g:X \to \Omega$ be the map $g(x)=\varphi^{-f(x)}(x)$. As a map $g:X \to X$ it is homotopic to the identity via $g_t(x)=\varphi^{-tf(x)}(x)$ and hence $g^*TX$ is isomorphic to $TX$ as a symplectic (or complex, or unitary) bundle. But $g^*TX$ is classified by the composition of the map $g:X \to \Omega$ and the map $\Omega \to BU(n)$ which classifies $T\Omega$. Hence $g^*TX$ admits an $\cI$-structure and therefore so does $TX$. \end{proof}

\subsection{Collaring co-isotropics}\label{sec: coll}
Before we study the singularities of the skeleta of Weinstein manifolds, it will be useful to record some results about embeddings of manifolds with corners in their cotangent bundles. Throughout this subsection, let $M$ denote an $n$-dimensional manifold with corners.

\begin{notation} We use the notation $\partial_kM$ to denote the locus of points modeled on 
\[0 \in [0,1)^k \times \bR^{n-k}.\] \end{notation}

\begin{definition}
Let $x \in \partial_k M$. We say that a 1-dimensional linear subspace $\ell \subset T_xM$ is a {\em collaring line} for $M$ at $x$ if there exists a vector $u \in \ell$ which is inwards pointing for $M$. 
\end{definition}

\begin{figure}[h]
\includegraphics[scale=0.07]{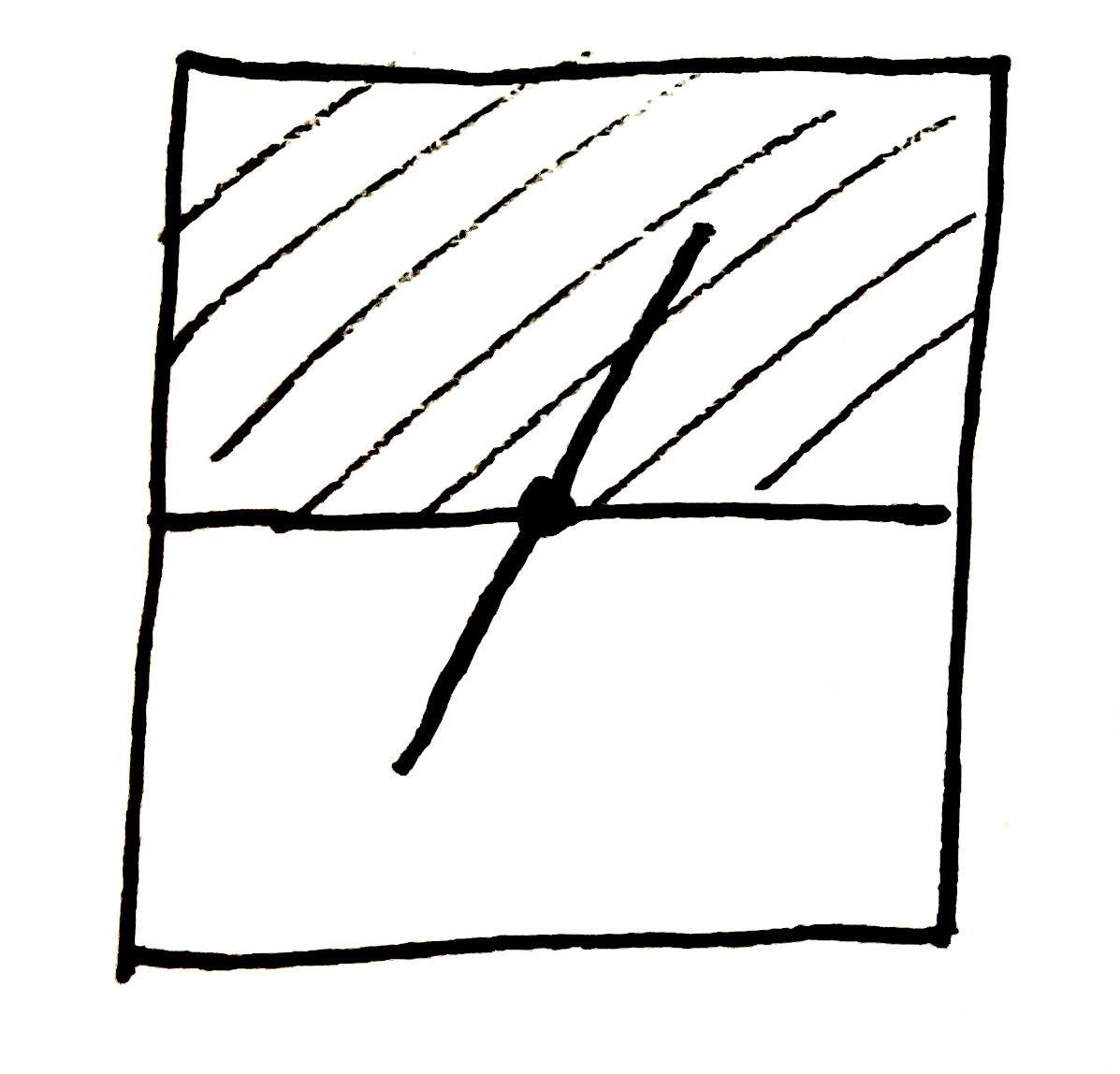}
\caption{A collaring line for a point at the boundary of a half-space.}
\label{fig:collar-1}
\end{figure}

\begin{figure}[h]
\includegraphics[scale=0.07]{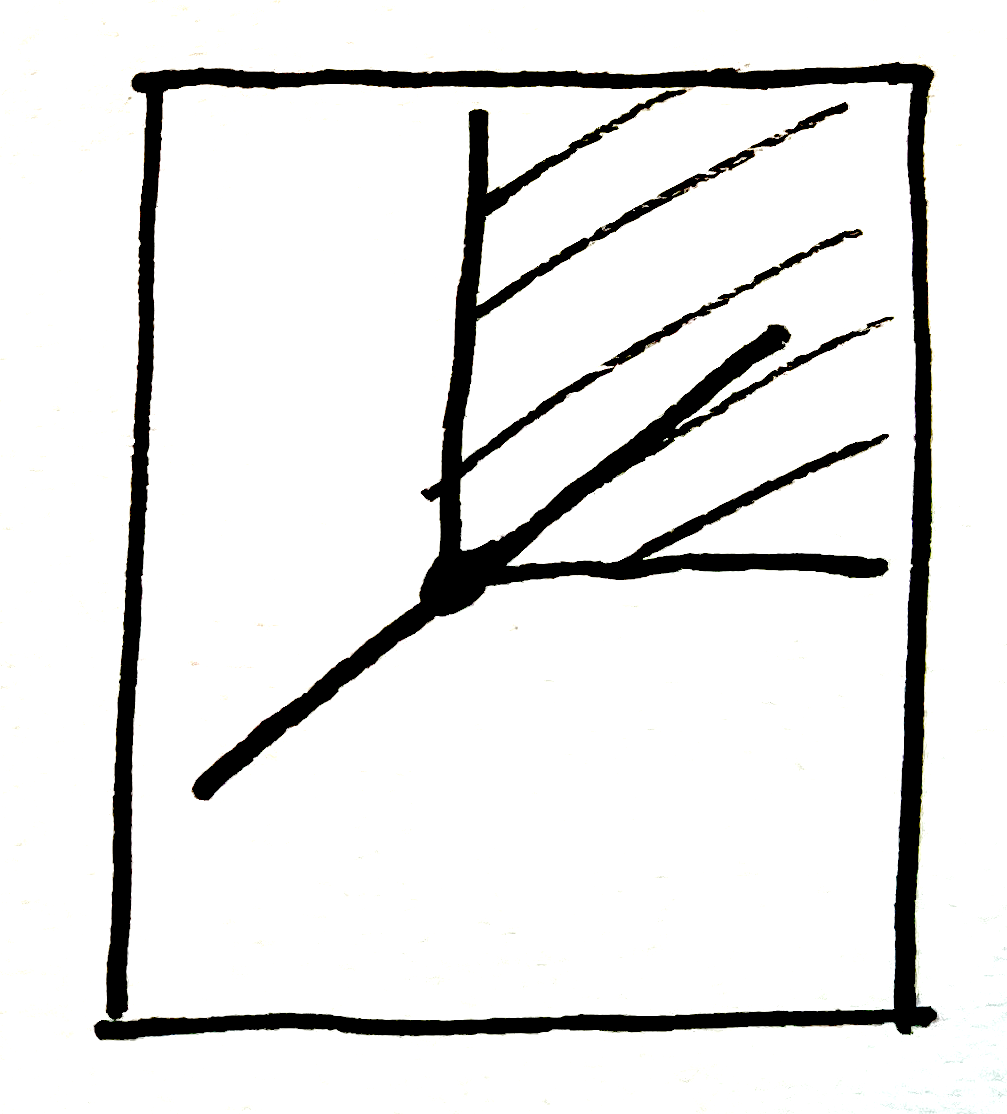}
\caption{A collaring line for a point at the boundary of the intersection of two half-spaces.}
\label{fig:collar-2}
\end{figure}

\begin{example}
For $x=0$ in $M= [0,1)^k \times \bR^{n-k}$, a line $\ell \subset \bR^n \simeq T_xM$ is collaring for $M$ at $x$ if and only if $\ell = \text{span}(u)$ for some non-zero vector $u =(u_1,\ldots,u_n) \in \bR^n$ such that $u_i>0$ for $i=1,\ldots, k$, see Figures \ref{fig:collar-1} and \ref{fig:collar-2}.
\end{example}

 \begin{definition} We say that a rank 1 distribution $\cE \subset TM|_{\partial M}$ is a  {\em collaring line field} if $\cE_x$ is a collaring line for $M$ at every $x \in \partial M$. 
 \end{definition}
 
 \begin{lemma}\label{lem:colar-line-contr} The space of collaring line fields on $\partial M$  is contractible. \end{lemma}
 
 \begin{proof} From the convexity of the condition $dx_i>0$ it readily follows that the space of collaring lines $\ell \subset T_xM$ for $M$ at a point $x \in \partial_kM$ is contractible. Further, if $\ell_x \subset T_xM$ is collaring at $x$, let $\ell$ be any continuous extension of $\ell_x$ to $TM|_U$ for an open subset $U \subset M$ containing $x$. Then the restriction of $\ell$ to a possibly smaller neighborhood $x \in V \subset U$ will be collaring at each $y \in V \cap \partial M$ (including points $y \in V \cap \partial_j M$ for $j \leq k$). From this the contractibility of any spheroid $\cE^z$, $z \in S^m$, of collaring line fields can be proved inductively using the stratification $\partial M = \coprod_{j=1}^k\partial_jM$, working from the deepest stratum $\partial_kM$ up.  
 \end{proof} 

\begin{definition} We say that an $(n+1)$-dimensional co-isotropic distribution $\cH \subset T(T^*M)|_{\partial M}$ is {\em collaring} if 
\begin{enumerate}
\item $\cH \pitchfork TM$, and
\item the intersection $\cE=\cH \cap TM \subset TM|_{\partial M}$ is a collaring line field for $M$. 
\end{enumerate}
When the collaring condition for a line field or a distribution holds along a subset $A \subset \partial M$, we will say that the line field or co-isotropic distribution is collaring over $A$.
\end{definition}

\begin{lemma}\label{lem:colar-contr-univ}
Let $(W,\omega)$ be a $2n$-dimensional symplectic vector space, $L \subset W$ a Lagrangian subspace and $\ell \subset L$ a line. The space $\cU=\cU_{L,\ell}(W)$ of $(n+1)$-dimensional co-isotropic subspaces $U \subset W$ such that $U \cap L = \ell$ is contractible.
\end{lemma}

\begin{proof} 
Given $U \in \cU$, first note that since $U \pitchfork L$ we have
$$(U^{\omega_\perp} \cap L)^{\omega_\perp} = (U^{\omega_\perp})^{\omega_\perp} + L^{\omega_\perp} = U+L = W $$
and therefore $U^{\omega_\perp} \cap L = 0$.  

Given a Lagrangian $K \subset W$ transverse to $L$ we get an $(n+1)$-dimensional co-isotropic $U=K + \ell$, indeed $U^{\omega_\perp} = K^{\omega_\perp} \cap \ell^{\omega_\perp} = K \cap \ell^{\omega_\perp} \subset K \subset U$. This defines a map from the contractible space $\Lambda_L(W)$ of Lagrangians $K \subset W$ transverse to $L$ to the space $\cU$. This map is readily verified to be a Serre fibration and therefore it suffices to check that the fiber is contractible.

For this purpose fix a co-isotropic $U \subset W$ such that $U \cap L= \ell$. If $K \subset U$ is a Lagrangian subspace transverse to $L$, then $U^{\omega_\perp} \subset K^{\omega_\perp}=K$ and $U=K+\ell$. Therefore, $\lambda = K/U^{\omega_\perp} \subset U/U^{\omega_\perp}$ is a line transverse to the line $(U^{\omega_\perp}+\ell)/U^{\omega_\perp}$ in $U/U^{\omega_\perp}$. Conversely, any such line corresponds to a unique Lagrangian $K \subset U$ transverse to $L$. The fiber thus becomes identified with the set of lines in a 2-dimensional space transverse to a fixed line, which is contractible. \end{proof}
 
 \begin{lemma}\label{lem:cois-contr-loc}
 Consider the submanifold with corners $M=[0,\infty)^k \times \bR^{n-k}$ where we identify the symplectic vector spaces $T_0(T^*M)$ with $\bC^n$ in the usual way $z_j =q_j +i p_j$. The space of $(n+1)$-dimensional co-isotropic subspaces $U \subset \bC^n$ which are collaring for $M$ at the origin is contractible. 
 \end{lemma}
 
 \begin{proof}
 The assignment $U \mapsto U \cap \bR^n$ defines a map from the space of $(n+1)$-dimensional co-isotropic subspaces $U \subset \bC^n$ which are collaring for $M$ at the origin to the space of collaring lines for $M$ at the origin. This map is a Serre fibration, the target is contractible and the fiber over a collaring line $\ell \subset \bR^n$ is $\cU_{\bR^n,\ell}(\bC^n)$, so we are done by Lemma \ref{lem:colar-contr-univ}.  \end{proof}

\begin{proposition}\label{prop:coll-contr}
Let $M$ be a manifold with corners. The space of collaring co-isotropic distributions on $\partial M$ is contractible.\end{proposition}

\begin{proof}
Once again the contractibility of any spheroid $\cH^z$, $z \in S^m$ of collaring co-isotropic distributions may be proved inductively using the stratification $\partial M = \coprod_{j=1}^k \partial_k M$. The necessary observation in this case is the fact that if $\cH_x \subset T_x(T^*M)$ is a collaring co-isotropic at $x \in \partial M$, and $\cH$ is any continuous extension to a co-isotropic distribution on an open subset $U \subset M$ containing $x$, then on a possibly smaller open subset $x \in V  \subset U$ the co-isotropic $\cH_y \subset T_y(T^*M)$ is collaring at all $y \in V \cap \partial M$. 
 \end{proof}

\subsection{Arboreal singularities}
\begin{figure}[h]
\includegraphics[scale=0.1]{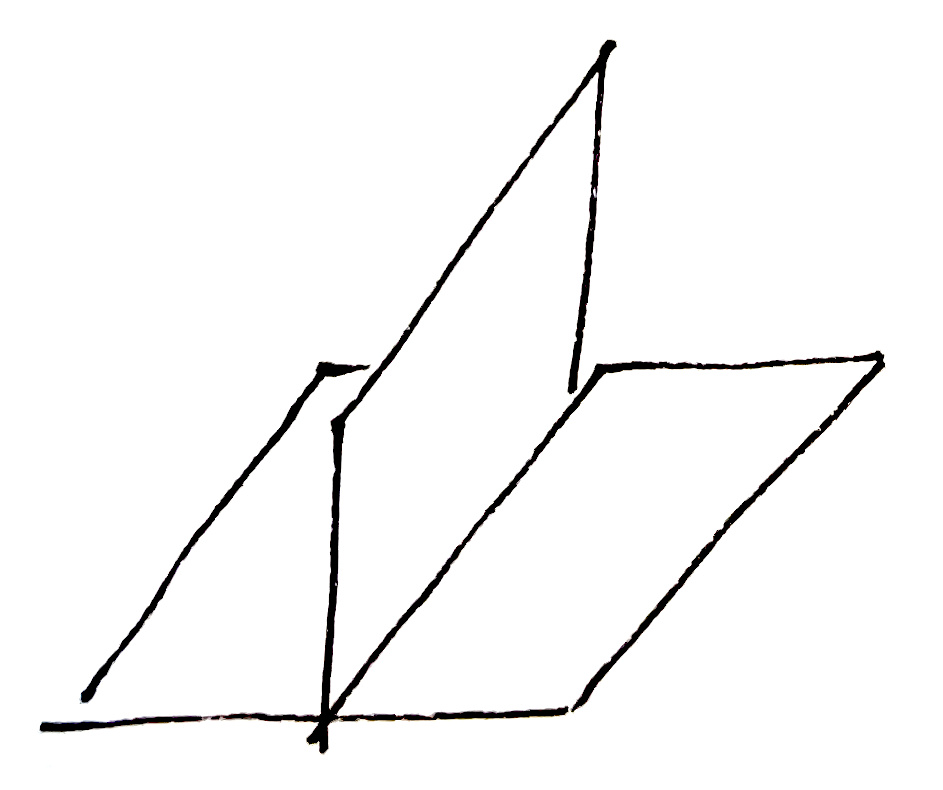}
\caption{The $A_2$ arboreal singularity.}
\label{fig:A_2}
\end{figure}

\begin{figure}[h]
\includegraphics[scale=0.2]{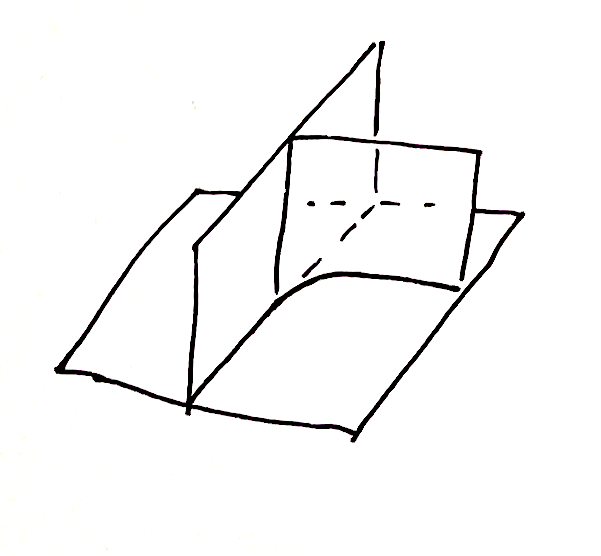}
\caption{The $A_3$ arboreal singularity.}
\label{fig:A_3}
\end{figure}

\begin{figure}[h]
\includegraphics[scale=0.08]{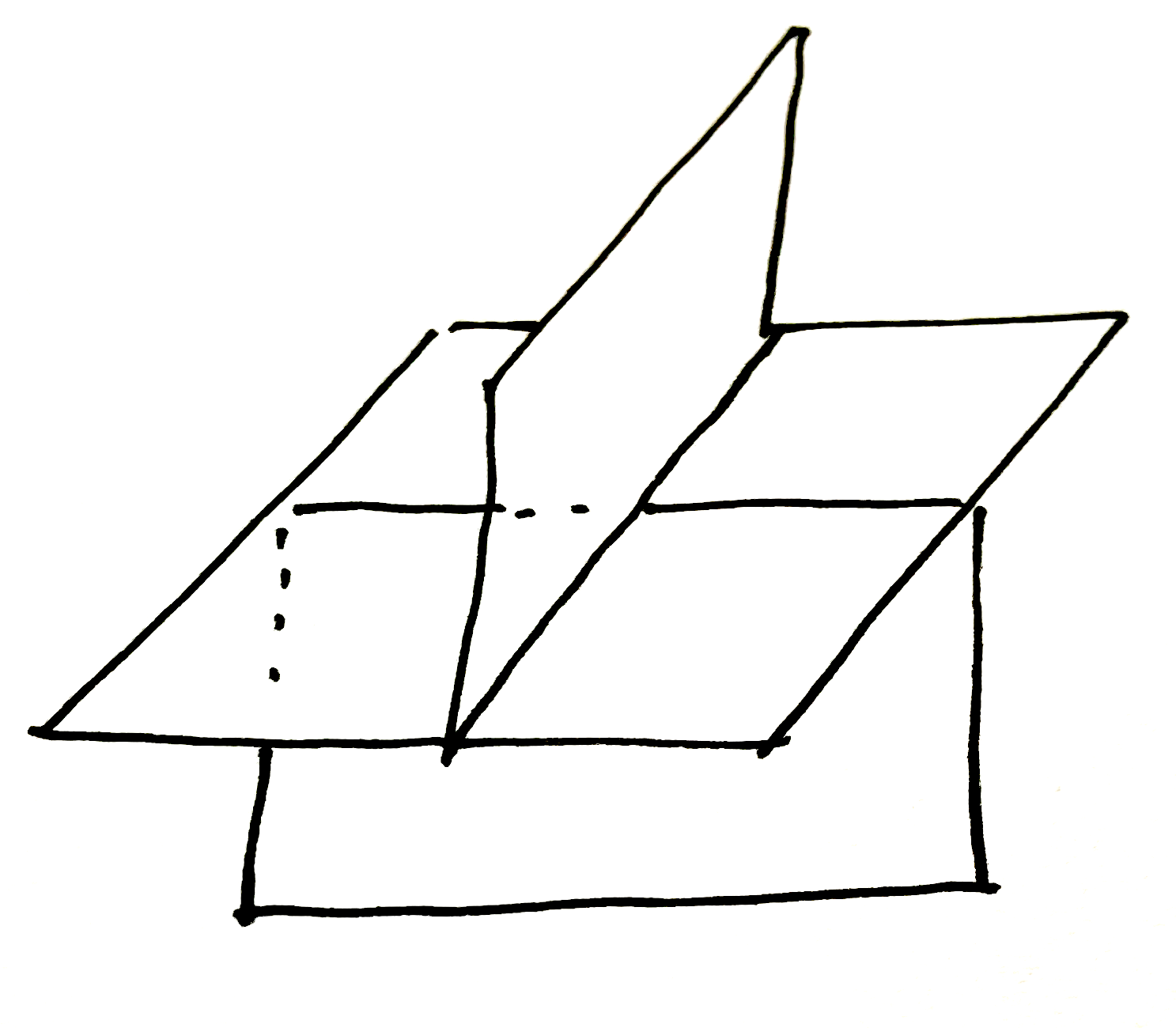}
\caption{The arboreal singularity corresponding to the $A_3$ tree but with the root at the central vertex.}
\label{fig:2D-arb-sing}
\end{figure}

In general, the singularities of the skeleton of a Weinstein manifold can be quite complicated. The class of arboreal singularities was introduced by Nadler in \cite{Na1} as a reasonable class of Lagrangian singularities for the skeleta of Weinstein manifolds, see Figures \ref{fig:A_2}, \ref{fig:A_3}, and \ref{fig:2D-arb-sing} which illustrate a few small dimensional examples. The class of arboreal Lagrangians is the smallest class of germs of singular Lagrangians characterized up to symplectomorphism by the following axioms, where we denote $\nu = \ker (d \pi)$ the vertical distribution of the fibration $\pi: S^*\bR^n \to \bR^n$.

\begin{enumerate}
\item A point is an arboreal Lagrangian in $T^*\bR^0=\bR^0$.
\item If $L \subset T^*\bR^n$ is an arboreal Lagrangian, then $L \times \bR \subset T^*\bR^n \times T^*\bR \simeq  T^*\bR^{n+1}$ is also an arboreal Lagrangian, where $L \times \bR$ denotes the product of $L$ with the zero section $\bR \subset T^*\bR$.
\item If  $\Lambda = \coprod_i \Lambda_i \subset S^*\bR^n$ is a Legendrian embedding of a disjoint union of arboreal Lagrangian germs at points $x_i \in \pi^{-1}(0)$ such that $\Lambda \pitchfork \nu$, and if moreover the projections $\Lambda_i \xrightarrow{\pi} \bR^n$ are mutually transverse (in the strongest possible sense), then the germ at the origin of $L=\bR^n \cup \text{Cone}(\Lambda) \subset T^*\bR^n$ is an arboreal Lagrangian.
\end{enumerate}

\begin{figure}[h]
\includegraphics[scale=0.2]{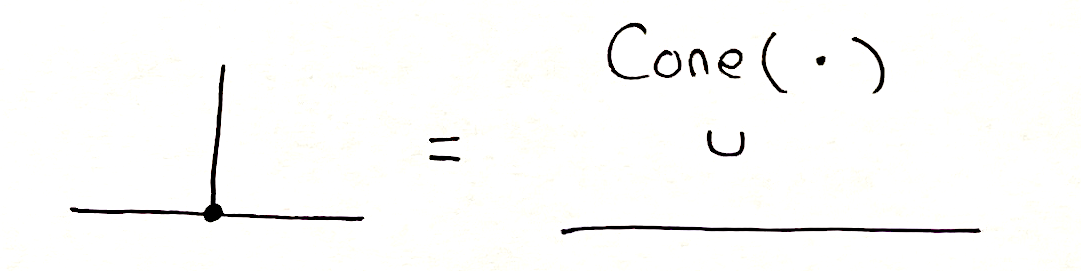}
\caption{The arboreal singularity of Figure \ref{fig:A_2} is the union of the fiberwise Liouville cone on a Legendrian $\Lambda \subset S^*\bR^n$ and the zero section $\bR^n$.}
\label{fig:A_2-cone}
\end{figure}

\begin{figure}[h]
\includegraphics[scale=0.2]{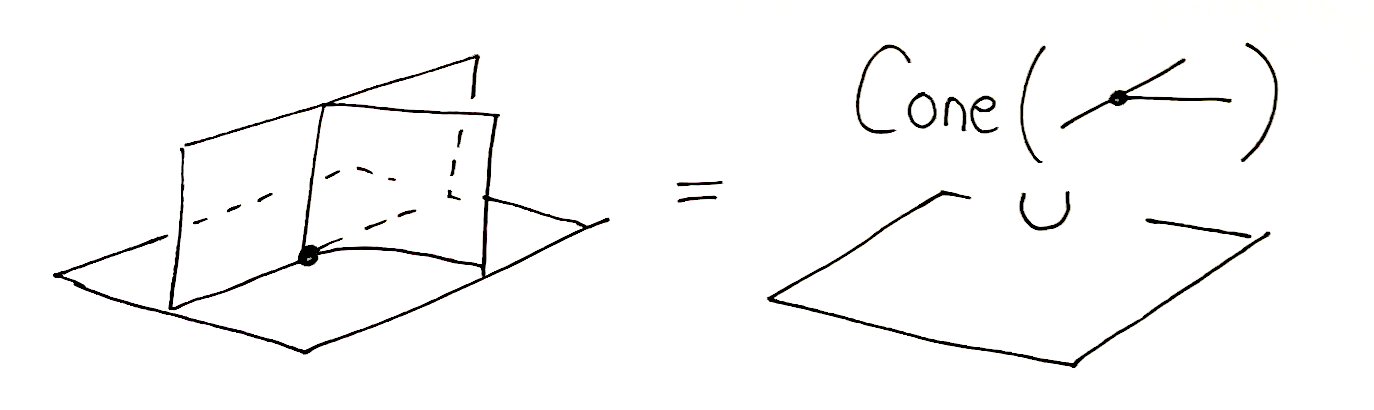}
\caption{The arboreal singularity of Figure \ref{fig:A_3} is the union of the fiberwise Liouville cone on a Legendrian $\Lambda \subset S^*\bR^n$ and the zero section $\bR^n$.}
\label{fig:A_3-cone}
\end{figure}

\begin{figure}[h]
\includegraphics[scale=0.2]{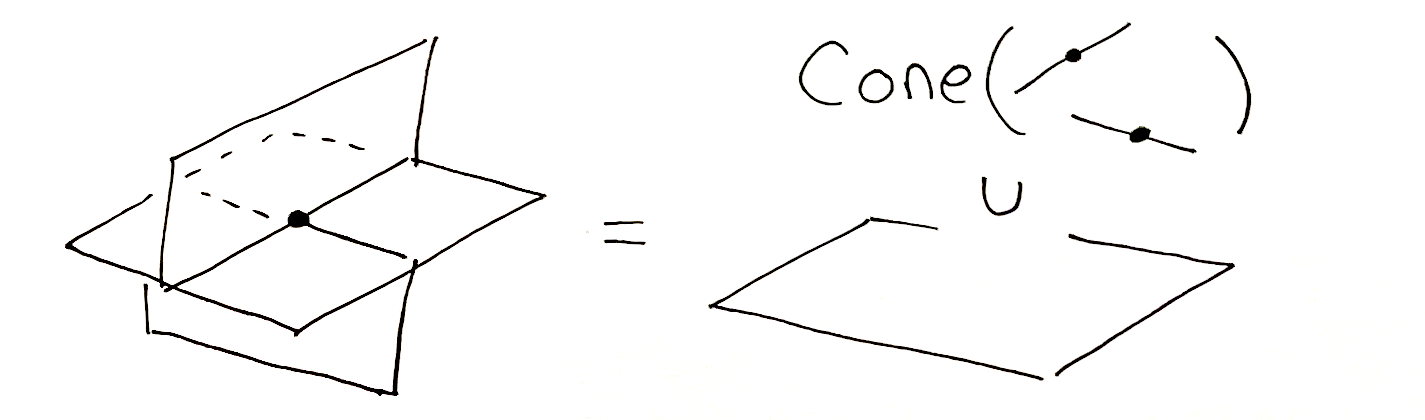}
\caption{The arboreal singularity of Figure \ref{fig:2D-arb-sing} is the union of the fiberwise Liouville cone on a Legendrian $\Lambda \subset S^*\bR^n$ and the zero section $\bR^n$.}
\label{fig:2D-cone}
\end{figure}

\begin{remark}
Here $\text{Cone}(\Lambda)$ denotes the Liouville cone of $\Lambda$ with respect to the canonical Liouville structure $pdq$, see Figures \ref{fig:A_2-cone}, \ref{fig:A_3-cone} and \ref{fig:2D-cone} for illustrations.
\end{remark}

To each germ of arboreal singularity is naturally associated a triple $\mathcal{T}=(T,\rho,\varepsilon)$ where $T$ is a finite connected tree, $\rho$ is a root for $T$, and $\varepsilon$ is a choice of sign $\pm 1$ for each edge not adjacent to the root. It is proved in \cite{AGEN1} that this combinatorial data uniquely characterizes the arboreal singularity germ up to symplectomorphism, and moreover, the space of symplectomorphisms preserving a given arboreal singularity is contractible. Thus  arboreal singularities are completely characterized by the combinatorial data $(T,p,\varepsilon)$. 

We say that a Weinstein manifold $(X,V,\phi)$ admits an arboreal skeleton if there is a homotopy of Weinstein structures $(V_t,\phi_t)$ such that $(V_0,\phi_0)=(V,\phi)$ and $(X,V_1,\phi_1)$ has an arboreal skeleton. In \cite{AGEN3} it is shown that polarizable Weinstein manifolds always admit arboreal skeleta. In fact there holds the following equivalence:

\begin{theorem}[\cite{AGEN3}] A Weinstein manifold $X$ admits a skeleton whose singularities are all symplectomorphic to  positive arboreal models $(T,p,\varepsilon)$, $\varepsilon_i=+1$, if and only if $X$ admits a polarization. \end{theorem}

There are currently no known sufficient conditions for the existence of an arboreal skeleton in high dimensional Weinstein manifolds beyond the positive case settled by the above result. The goal of this section is to establish a necessary condition, namely the existence of an $\cI$-structure on the tangent bundle of any arborealizable $X$. 

\subsection{Geometry of arboreal singularities}\label{sec: geometry}

In this section we derive some basic properties of arboreal singularities. We rely throughout on the explicit arboreal models of \cite{AGEN1}.

\begin{definition} Any arboreal singularity is of the form $K=\bR^n \cup \text{Cone}(\Lambda) \subset T^*\bR^n$ where $\bR^n \subset T^*\bR^n$ is the zero section and $\Lambda \subset S^*\bR^n$ is an arboreal Legendrian. The germ of the smooth Lagrangian submanifold $\bR^n \subset T^*\bR^n$ at the origin will be called the {\em local zero section} of the arboreal singularity $K$.
\end{definition}

\begin{definition}
Let $\cK \subset W$ be an arboreal skeleton of a Weinstein manifold $W$. A {\em smooth piece} of $\cK$ is a maximal smooth Lagrangian submanifold $L \subset \cK$ which is given by the union of the local zero sections of the arboreal singularities it encounters.
\end{definition}

Smooth pieces will not in general be relatively closed or compact. They are disjoint and their union is all of $\cK$.

\begin{example} The local model for the $A_2$ singularity within $T^*\bR$ is composed of two smooth pieces: the zero section, i.e., $\bR \subset T^*\bR$; and the set $\{(0,p) \mid p \in \bR, p > 0\}$. \end{example}

 Given a smooth piece $L$ of $\cK$, there is a unique compactification of $L$ to a manifold with corners $L^c$ such that there is a canonical diffeomorphism $L^c \setminus \partial L^c \to L$ and such that the embedding $L \subset X$ extends to an immersion $\iota_L: L^c \to X$ with image $\overline{L}$. We call $L^c$ the {\em completion} of $L$. Note that $L^c \to X$ need not be an embedding. We will sometimes refer to a point $p \in L$ as being in the interior of the smooth piece $L$, even if this is redundant, for emphasis.

Now fix a signed rooted tree $\mathcal{T}=(T,\rho,\ve)$. We let  $v(T)$ denote the set of vertices of $T$, and define $a(T) \subset v(T)$ to be the subset of vertices adjacent to the root $\rho$. After enumerating the elements of $a(T)$ by $v_1, \ldots, v_p$, let $\mathcal{T}_i =(T_i, v_i, \ve_i)$, where $T_i$ is the subgraph obtained by removing the connected components of $T \setminus \{\rho\}$ which do not contain $v_i$, and $\varepsilon_i$ is obtained by restricting $\varepsilon$ to $T_i$. We  partition the non-root vertices of $T$ by writing $ v(T) \setminus \{\rho \}= \coprod_{i=1}^p v(T_i)$. 
 Let $i: v(T) \setminus \{\rho\} \to \{1, \ldots, p\}$ be defined according to this partition. Denote $n_i=\# v(T_i)$, the number of vertices of $T_i$.
 
Let $0 \in K =  \bR^n \cup \text{Cone}(\Lambda) \subset T^*\bR^n$ be the germ at the origin of the arboreal singularity corresponding to $\mathcal{T}.$ The number of vertices in $T$ is at most $n+1$, but for our analysis it is enough to consider the case where it is exactly $n+1$, as the other arboreal singularities are given by taking trivial products $K \times \bR^k \subset T^*\bR^n \times T^*\bR^k = T^*\bR^{n+k}$, where $\bR^k \subset T^*\bR^k$ is the zero section. We have $\Lambda = \bigcup_{i=1}^p \Lambda_i$ for arboreal Legendrian germs at distinct points $x_i \in \pi^{-1}(0) \subset  S^*\bR^n$, with $\Lambda_i \subset S^*\bR^n$ a (stabilization) of a Legendrian lift of the arboreal singularity $K_i \subset T^*\bR^{n_i}$ corresponding to the signed rooted tree $\mathcal{T}_i$. 

\begin{figure}[h]
\includegraphics[scale=0.1]{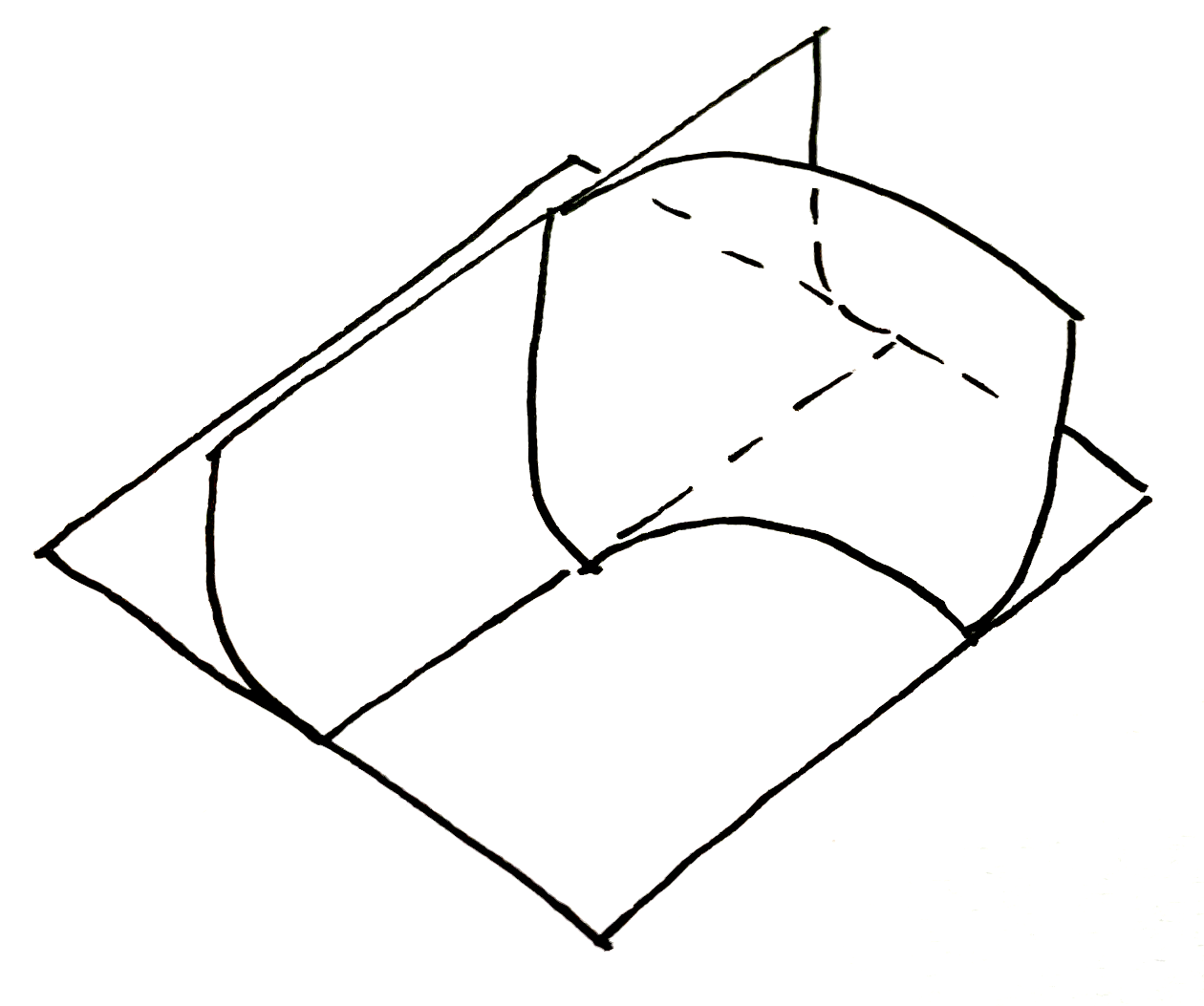}
\caption{The $A_4$ arboreal singularity in $T^*\bR^3$ is the conormal to the above (co-oriented) hypersurface of $\bR^3$.}
\label{fig:A_4-arb-hyp}
\end{figure}

Concretely, the Legendrian $\Lambda$ is the positive conormal of a singular co-oriented hypersurface $ \cA \subset \bR^n$, and according to the standard model of \cite{AGEN1}, up to symplectomorphism we can write $\cA = \bigcup_{i=1}^p \wh \cA_i$ where  $$\wh \cA_i = \bR^{n_1 } \times \cdots  \times \bR^{n_{i-1}} \times \cA_i \times \bR^{n_{i+1} } \times \cdots \bR^{n_p } \subset \bR^n$$ where $\cA_i \subset \bR^{n_i}$ is the front projection of the arboreal singularity $K_i$ corresponding to $T_i$ under the Legendrian embedding in $S^*\bR^{n_i}$ given by the product Legendrian $\Lambda_i$. Each $\wh \cA_i$ is piecewise smooth and has a well-defined co-oriented tangent plane field. We denote by $Q_i \subset T_0 \bR^n$ its tangent plane at the origin, which is the hyperplane $\{x_i=0\}$. The hyperplanes $Q_1,\ldots,Q_p \subset T_0\bR^n$ are mutually transverse in the strongest possible sense: each $Q_i$ is transverse to all possible intersections of the $Q_j$, $j \neq i$. Equivalently, the corresponding co-orientations $x_i \in S^*\bR^n$ are linearly independent as covectors in $T^*\bR^n$. See Figures \ref{fig:A_4-arb-hyp},  for an illustration.

\begin{figure}[h]
\includegraphics[scale=0.15]{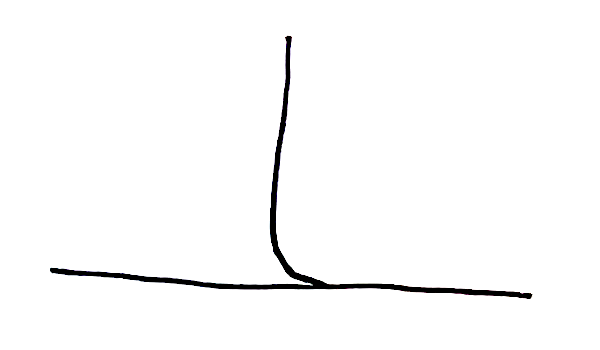}
\caption{The $A_3$ arboreal singularity in $T^*\bR^2$ is the conormal to the above (co-oriented) hypersurface of $\bR^2$.}
\label{fig:A_2-hyp}
\end{figure}

\begin{figure}[h]
\includegraphics[scale=0.13]{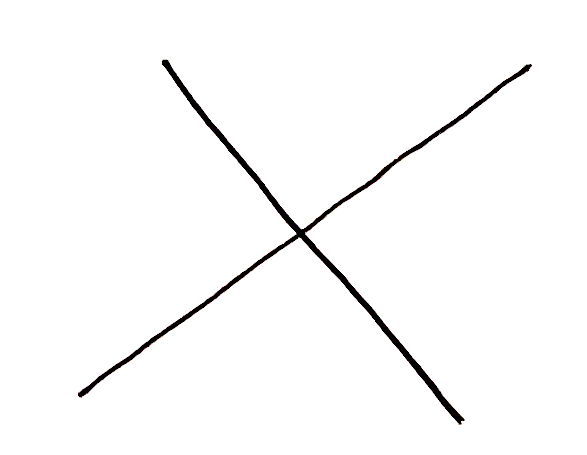}
\caption{The arboreal singularity in $T^*\bR^2$ of Figure \ref{fig:2D-arb-sing} is the conormal to the above (co-oriented) singular hypersurface of $\bR^2$. Note that the hypersurface is of the form $(\cA \times \bR) \cup (\bR \times \bA)$ for $\A=\{0\} \subset \bR$ the point hypersurface corresponding to the $A_2$ arboreal singularity in $T^*\bR$..}
\label{fig:2D-hyp}
\end{figure}

\begin{figure}[h]
\includegraphics[scale=0.08]{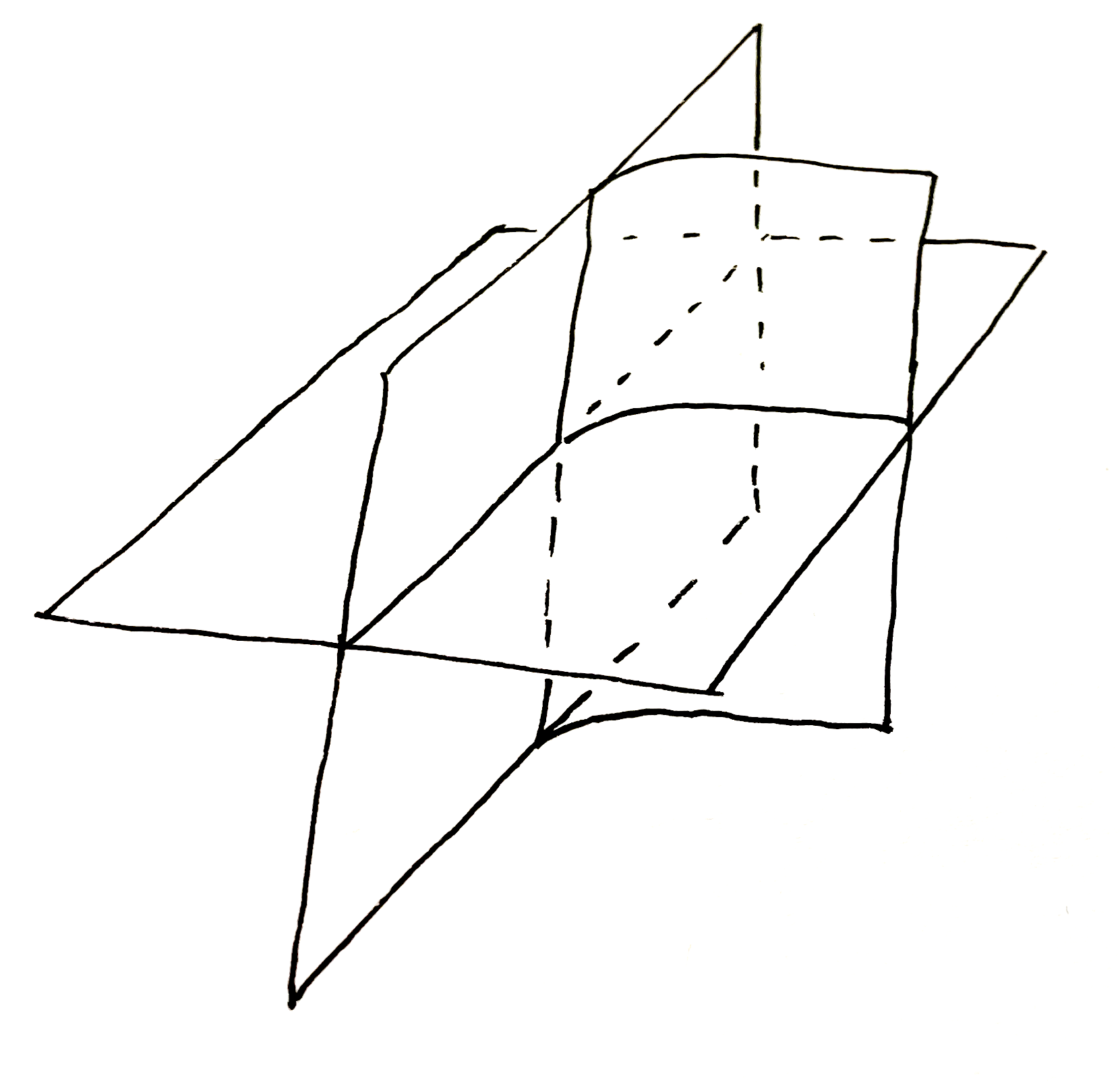}
\caption{Another arboreal singularity in $T^*\bR^3$ is given by the conormal to the above (co-oriented) hypersurface of $\bR^3$. Note that the hypersurface is of the form $(\cA \times \bR^2) \cup (\bR^2 \times \cA')$ for $\cA=\{0\} \subset \bR$ the hypersurface corresponding to the 1-dimensional $A_2$ singularity in $T^*\bR$ and $\cA' \subset \bR^2$ is the hypersurface of the $A_3$ singularity in $T^*\bR^2$ as in Figure \ref{fig:A_2-hyp}.}
\label{fig:3D-arb-hyp}
\end{figure}

For each $i$, the smooth pieces that make up  $\Lambda_i$ are in bijection with the vertices of $T_i$; for $v \in v(T_i)$ we write $L_v$ for the corresponding smooth piece.  These give rise to smooth pieces $\wh L_{v} = \text{Cone}(L_{v} )\subset S^*\bR^n \setminus \bR^n$ for the arboreal Lagrangian $L$. Note that we have a diffeomorphism $\wh L_{v} \simeq L_{v} \times [0,1)$ where $[0,1)$ parametrizes the Liouville direction. Each smooth piece $\wh L_{v}$ is the positive conormal to a smooth co-oriented hypersurface in $\bR^n$ (the union of which over all $v \in v(T_i)$ make up $\wh \cA_i$). Further note that $\wh L_{v}^c \simeq L_{v}^c \times [0,1]$ and $\iota_{\wh L_{v} } :   L_{v}^c \times [0,1] \to T^*\bR^n$ sends a point $(y,1) \in L_{v} \times \{ 1\}$ to $\pi( \iota_{L_{v} }(y)) \in  \wh \cA_{i(v)} \subset \bR^n$. For each $v \in v(T) \setminus \{\rho\}$, denote   
\[\rhat{L_v} = \iota_{\wh L_{v} } (L_{v}^c \times [0,1]). \]
This is an embedded manifold with corners. 

When two vertices $v$ and $v'$ of $T$ lie within the same subgraph $T_i$, both $\rhat{L_v}$ and $\rhat{L_{v'}}$ intersect $\pi^{-1}(0)$ at the same point $x_i \in S^*\bR^n$; thus the tangent spaces $T_0(\, \rhat{L_v}\, )$ and $T_0(\,\rhat{L_{v'}}\,)$, regarded as subspaces of $T_0(T^* \bR^n)$, coincide. We can thereby define a collection of subspaces 
\[ P_i \subset T_0(T^*\bR^n),  i=1, \ldots, p.\]   Explicitly, with respect to the isomorphism  $T_0(T^*\bR^n) \simeq T_0\bR^n \oplus T^*_0\bR^n$, $P_i$ is the direct sum $Q_i \oplus \ell_i$ of the hyperplane $Q_i \subset T_0\bR^n$ and the line $\ell_i$ in $T_0^*\bR^n$ of covectors at the origin which are zero on $Q_i$, which is spanned by the co-orientation $x_i$ of $Q_i$.

Note that for each $1 \leq i \leq p$ indexing a root-adjacent vertex $v_i \in a(T)$,  the piece $L_{v_i}$ is the zero section for the arboreal singularity at $x_i \in \Lambda_i$, and $\rhat{L_{v_i}}$ is a manifold with boundary, whose boundary is contained within the zero section $\bR^n \subset T^* \bR^n$: 
\[ \del \rhat{L_{v_i}} = \iota_{L_{v_i}}(L_{v_i} \times \{1\}) = \{x_i = 0\} \subset \bR^n.\]
Thus the space of vectors in $P_i$ which are inwards pointing for $\rhat{L_{v_i}}$ at the origin consists of the half space $P_i^+=Q_i \oplus \ell_i^+$ where $\ell_i^+ \subset \ell_i$ is the set of covectors at the origin which are zero on $Q_i$ and endow $Q_i$ with the same co-orientation as $x_i$, i.e., positive multiples of $x_i$. More generally we have the following lemma:

 \begin{lemma} \label{lem:half-space}
For fixed $1 \leq i \leq p$, the space of vectors in $P_i \subset T_0(T^*\bR^n)$ which are simultaneously inwards pointing for all of the submanifolds with boundary $\rhat{L_v}$ as  $v$ ranges over the vertices of $T_i$ is an intersection of half-planes 
\[ \bigcap_{k=1}^{n_i} \{ \nu_{i}^k > 0 \} \subset P_i\] where $Q_i = \{ \nu_i^1 = 0 \} \subset P_i$ and $\nu_i^2,\ldots,\nu_i^{n_i}$ are a set of linearly independent elements of $P_{i}^*$ which pull back to a set of linearly independent elements of $Q_i^*$. In particular, $\nu_i^1,\nu_i^2,\ldots,\nu_i^{n_i}$ are linearly independent in $P_i^*$ and so this space is contractible.
\end{lemma}

\begin{proof}
The proof is by induction on the dimension of the arboreal singularity. If we assume that the conclusion holds for the arboreal singularity $\mathcal{T}=(T,\rho,\ve)$, it will suffice to show that it holds for the arboreal singularity $\mathcal{T}'=(T',\rho',\ve')$, obtained by adding a new vertex $\rho'$ to $T$, which will be the new root, an edge between the new vertex and the root, and an arbitrary extension of the signs $\ve$ to the edges adjacent to $\rho$ but not to $\rho'$. 

More geometrically, let us assume that the conclusion holds for $0 \in K= \bR^n \cup \text{Cone}(\Lambda)$, the arboreal singularity corresponding to $\mathcal{T}$.
Recall that $\wh \cA_i \subset \bR^n$ is a product $\bR^{n_1+\cdots + n_{i-1}} \times \cA_i \times \bR^{n_{i+1} + \cdots + n_p }$ and hence the covectors $\nu_i^k$ are zero on the subspace $\bR^{n_1+\cdots + n_{i-1} } \times 0^{n_i} \times \bR^{n_{i+1} + \cdots + n_p} \subset Q_i$. 
The standard model for the arboreal singularity for $\mathcal{T}'$ is obtained by taking a suitable Legendrian embedding $K \hookrightarrow S^*\bR^{n+1}$, $0 \mapsto \wh x \in \pi^{-1}(0)$ and forming an arboreal singularity $\bR^{n+1} \cup \text{Cone}(K)$ at $0$. The smooth pieces consist of the zero section $\bR^{n+1}$ together with the product of smooth pieces for $K$ with the new Liouville direction $[0,1)$, i.e., in terms of our previous notation they are of the form $\wh L_{v} \times [0,1) = L_{v} \times [0,1) \times [0,1)$. There is also the smooth piece $\bR^n \times [0,1)$ where $\bR^n$ is the zero section of $K=\bR^n \cup \text{Cone}(\Lambda)$, and $[0,1)$ is the new Liouville direction.

At the point $0 \in \bR^{n+1}$ the common tangent plane to the smooth pieces which are not the zero section $\bR^{n+1}$ of $\bR^{n+1} \cup \text{Cone}(K)$ are all equal to the conormal $P$ of the hyperplane $Q=\{ \wh x = 0 \} \subset T_0 \bR^{n+1}$. We may write $P \simeq Q \times \bR$ where $0 \times \bR^+$  is the ray of positive multiples of $\wh x$. In particular, quotienting out by the Liouville direction we obtain an identification of the tangent planes $P_i$ of the smooth pieces $\wh L_{v_i}$ at $0 \in T^*\bR^n$ with $Q$, with each other and with the tangent plane $T_0\bR^n$ of the zero section $\bR^n$ of $K$.  The identification $P_i \simeq T_0\bR^n$ fixes the hyperplane $Q_i$ which is common to both and under this identification the covector $\nu_i^1 \in P_i^*$ becomes the covector $\nu_i^1=x_i \in T_0^*\bR^n$ which still satisfies $Q_i = \{\nu_i^1=0\}$, and the other covectors $\nu_i^k \in T_0^*\bR^n$, $k>1$, still pull back to linearly independent elements of $Q_i^*$.

Now, the subspace of $P$ consisting of vectors which are inwards pointing to the smooth pieces $\wh L_{v} \times [0,1)$ and to the smooth piece $\bR^n \times [0,1)$ correspond respectively to the following two subsets of $P \simeq T_0\bR^n \times \bR$.
\begin{enumerate}
\item The subset of vectors in $P_{i(v)} \simeq T_0\bR^n$ which are inwards pointing to $\wh L_{v}$, multiplied by $\bR^+$, and
\item  the half-space $T_0\bR^n \times \bR^+$ in $P$ consisting of vectors which co-orient the hyperplane $Q=\{ \wh x = 0 \}$ positively, i.e., with the same co-orientation as $\wh x$.
\end{enumerate}
The conclusion then follows by applying the inductive hypothesis. Indeed, the above collection of covectors $\nu^i_k$ for all $i$ and $k$ remains linearly independent when put together in $T_0\bR^n$. To see this it suffices to observe that each subcollection $\nu_i^1,\nu_i^2,\cdots,\nu_i^{n_i} \in T_0^*\bR^n$ for fixed $i$ is zero on the subspace $\bR^{n_1+\cdots + n_{i-1} } \times 0^{n_i} \times \bR^{n_{i+1} + \cdots + n_p} \subset Q_i$ and linearly independent when restricted to $0^{n_1+\cdots+n_{i-1} } \times \bR^{n_i} \times 0^{n_{i+1} + \cdots + n_p} $. \end{proof}

\begin{remark} While we work with $n = \# v(T) -1$ above, the above proof can be easily modified to prove the same statement in the stabilized case  $n >  \# v(T)-1$. \end{remark}

 \subsection{Germs of collaring co-isotropics for arboreal singularities}
 
 We continue using the notation of Section \ref{sec: geometry}: we fix a triple $\mathcal{T}= (T, \rho, \epsilon)$ indexing a germ of an arboreal singularity $K \subset T^* \bR^n$ (here we allow the case where $n > \# v(T)-1$) and enumerate the elements of $a(T)$, i.e., the set of vertices adjacent to the root $\rho$, by $v_1, \ldots, v_p$. This defines a partitioning function $i: v(T) \setminus \{\rho\} \to \{1, \ldots, p \}$.  For each $v \in v(T) \setminus \{\rho\}$, we have a manifold with corners $\rhat{L_v}$ embedded within $T^* \bR^n$, and for each $1 \leq i \leq p$, we have a subspace $P_i \subset T_0 T^* \bR^n$ such that $T_0 \rhat{L_v}= P_{i(v)}$ for all $v \in v(T) \setminus \rho$.

From Lemma \ref{lem:half-space} the space of vectors in $P_i$ which point inwards to all $ \rhat{L_v}$, $v \in v(T_i)$, consists of a convex subset, and hence the space of lines in $P_i$ which are collaring for all $ \rhat{L_v}$, $j=1,\ldots, n_i$, is contractible. Arguing as in the previous section we conclude that the space of $(n+1)$-dimensional co-isotropic subspaces of $T_0(T^*\bR^n)$ which are collaring for all $\rhat{L_v}$, $v \in v(T_i)$, is also contractible. 

\begin{definition} 
For fixed $i=1,\ldots, p$, let $\bH_i$ be the space of germs at the origin of $(n+1)$-dimensional co-isotropic distributions on $T^*\bR^n$ which are collaring for all $\rhat{L_v}$, $v \in v(T_i)$.  The space of {\em collaring co-isotropic germs} for the arboreal singularity $K \subset T^* \bR^n$ is denoted
\[ \bH^{\mathcal{T},n} = \prod_{i=1}^p \bH_i,\] so a collaring co-isotropic germ consists of a tuple $(H_1,\ldots,H_{p})$ where $H_i \in \bH_i$.   \end{definition}

\begin{remark}\label{rem: aut trans} Let $(H_i)_{i} \in \bH^{\mathcal{T},n}$. We know that the intersection of $H_{i}$  and the Lagrangian plane $ T(\, \rhat{L_{v_i}}\,)$ at 0 is a line transverse to $Q_i$, which implies that $\mathrm{dim} (H_{i})_0 \cap T_{0} \bR^n=1$, i.e., $(H_i)_0$ is transverse to $T_{0} \bR^n$. By continuity, given any extension $\cF$ of $T_0\bR^n$ to a Lagrangian distribution on $T^*\bR^n$ we have that $H_i \pitchfork \cF$ on $Op(0) \cap \text{Cone}(\Lambda)$. 
\end{remark}

\begin{definition} \label{def: germs} We use the above local models to define the space of collaring co-isotropic germs $\bH_{x}$ at any point $x \in \cK$ of an arboreal Lagrangian $\cK$ in a symplectic manifold $X$, so if a point $x \in K$ is locally modeled on the singularity in $T^*\bR^n$ corresponding to a signed rooted tree $\mathcal{T}$, there is an identification $\bH_x \simeq \bH^{\mathcal{T},n}$. \end{definition}

We now want to make sense of how collaring co-isotropic germs can vary in families; to do this, we again specialize to the local model given above associated to the signed rooted tree $\mathcal{T}$ indexing the germ of a singularity $K \subset T^* \bR^n$ at the origin. Assume that $(H_i)_{i} \in  \bH^{\mathcal{T},n}$ is represented on an open neighborhood $U$ containing the origin. Let $y \in K \cap U$, and consider the signed rooted tree $(T_y, v_{i_y}, \varepsilon_y)$ indexing the local arboreal singularity at $y$. The vertices of $T_y$ can be regarded as a subset of the vertices of $T$: they index the set of $\rhat{L_v}$ which include $y$. We thus associate a set of isotropic distributions indexed by elements $w_1, \ldots, w_{p_y} \in a(T_y)$ by defining $H_{w_j} = H_{i(w_j)}.$ 

\begin{lemma}\label{lem:germ-ext}If $((H_i)_{i}, U)$ is a $p$-tuple of isotropic distributions on $U$ which defines an element of $\bH^{\mathcal{T},n}$, then exists some smaller neighborhood $V \subset U$ such that $0 \in V$ and $(H_{i})_i$ defines an element of $\bH_{y}$ for all $y \in V$ via the restriction procedure described above. \end{lemma}
\begin{proof}This follows from the openness of the collaring condition: we know that for any vertex $v \in v(T) \setminus \{0\}$, $H_{i(v)}$ is a collaring isotropic for $\rhat{L_{v}}$ at the origin, and so $H_{i(v)}$ remains collaring for $\rhat{L_{v}}$ in some neighborhood $U_v$ about the origin. We can then take $V = \bigcap_{v \in v(T) \setminus \{\rho\}} U_v$. \end{proof}

\begin{definition}\label{def: fam loc}
For an arboreal singularity $0 \in K = \bR^n \cup \text{Cone}(\Lambda) \subset T^*\bR^n$, a {\em holonomic family of collaring co-isotropic germs} defined in a neighborhood $U \subset K$ of the origin consists of the choice of a collaring co-isotropic germ at each $y \in U$ which arises by a single choice of a $p$-tuple of co-isotropic distributions defined on $U$ as in the conclusion of Lemma \ref{lem:germ-ext}. 
\end{definition}

\begin{definition}\label{def: fam glob}
For an arboreal Lagrangian $\cK$ in a symplectic manifold $W$, a {\em holonomic family of collaring co-isotropic germs} over $\cK$ consists of the choice of a collaring co-isotropic germ for each arboreal singularity $K_x=K \cap \Op(x)$, $x \in K$, which is a holonomic family of collaring co-isotropic germs near each point in the sense of Definition \ref{def: fam loc}.
\end{definition} 

\begin{remark} In general, one can define a sheaf $\bH_{\cK}$ supported on an arboreal Lagrangian $\cK$ by setting $\bH_{\cK}(U)$ to be the set of all holonomic families of collaring co-isotropic germs over $U$. Lemma \ref{lem:germ-ext} states that the stalk of $\bH_{\cK}$ at a point $x$ agrees with the space $\bH_{x}$. \end{remark}

\begin{lemma}\label{lem:hol-exist}
Let $\cK$ be an arboreal Lagrangian in a symplectic manifold $W$. Holonomic families of collaring co-isotropic germs over $\cK$ exist and are unique up to contractible choice of homotopy through such. 
\end{lemma}

\begin{proof} Existence may be proved as follows. Start with the deepest (top codimensional) stratum in the stratification $\cK = \cK^0 \supset \cK^1 \supset \cdots \supset \cK^n$ by codimension of the singularity type, where $\dim W = 2n$. For $x \in \cK^n$ a point, the space of collaring co-isotropic germs at $x$ is contractible, so we may choose such a collaring co-isotropic germ arbitrarily (the holonomic condition is vacuous at a point). 

Inductively, assume one has constructed a holonomic family of collaring co-isotropics over $\cK^{j+1}$ and let $R$ be an $(n-j)$-dimensional component of $\cK^j \setminus \cK^{j+1}$. Note that $R$ consists of arboreal singularities of the same combinatorial type. By Lemma \ref{lem:germ-ext}, the already defined family extends to a holonomic family on a  neighborhood of $\cK^{j+1}$ in $\cK^j$, so in particular we may assume we have a holonomic family of co-isotropic germs defined near $\partial R$. One may then extend this family to the rest of $R$ using the contractibility of the space of collaring co-isotropics for the fixed arboreal singularity type common to $R$. 

By proceeding in this way over each component $R$ we successfully extend the section to $\cK^j$, and the induction step is complete. The uniqueness up to contractible choice of homotopy may be proved in the same way. \end{proof}

 \begin{definition}\label{def: bf to em}
Given a collaring co-isotropic germ $H=(H_i)_i \in \bH$ at an arboreal singularity $0 \in K= \bR^n \cup \text{Cone}(\Lambda)$, represented on some open subset $U \subset T^*\bR^n$ containing the origin, we define an {\em associated} co-isotropic distribution $\cH$ defined on $V=U \cap (K \setminus \bR^n ) =  U \cap \text{Cone}(\Lambda)$ by setting $\cH=H_i$ on $U \cap \text{Cone}(\Lambda_i)$. Note that $\cH$ is not defined at the origin. \end{definition}

\begin{remark}
Note that $\cH$ is transverse to each smooth piece $L \subset \text{Cone}(\Lambda)$ of $K$ by the collaring condition, possibly after shrinking $V$ if necessary. Note also that if $\cF$ is any extension of $T\bR^n$ to a Lagrangian distribution on a neighborhood  of the zero section, then $\cF \pitchfork \cH$ on $V$, possibly after shrinking $V$ further if necessary, see Remark \ref{rem: aut trans}.
\end{remark}

\begin{remark}\label{rem: associated} Given a holonomic family of collaring co-isotropic germs for $\cK \subset W$ an arboreal Lagrangian in a symplectic manifold, the associated distributions $\cH^x$ are defined on relatively open subsets $V_x \subset K$ and by definition of the holonomic family condition we may assume that $\cH^x=\cH^y$ on $V_x \cap V_y$ by shrinking the $V_x$ if necessary. So we get an associated co-isotropic distribution $\cH$ defined on a relatively open subset of the form $V=\bigcup_x V_x \subset \cK$. Note that $\cK \setminus V$ is a disjoint union $\coprod_L L_V$ where $L$ ranges over the smooth pieces of $K$ and $L_V = L \setminus V$ is the complement in $L^c$ of a neighborhood of $\partial L^c$. \end{remark}

\subsection{Existence of $\cI$-structures for arboreal skeleta}

The main goal of this section is the following result:

\begin{theorem}\label{thm:I-str}
If a Weinstein manifold $X$ admits an arboreal skeleton, then the tangent bundle $TX$ admits an $\cI$-structure. 
\end{theorem}

\begin{proof}

Invoke Lemma \ref{lem:hol-exist} to produce a holonomic family of collaring co-isotropic germs along $\cK$ and let $\cH$ be associated the co-isotropic distribution  as in Remark \ref{rem: associated}. Recall that $\cH$ is defined   on a relatively open subset of the form $V=\bigcup_x V_x \subset \cK$, with $\cK \setminus V$ a disjoint union $\coprod_L L_V$ where $L$ ranges over the smooth pieces of $K$ and $L_V=L \setminus V$ is the complement in $L^c$ of a neighborhood of $\partial L^c$. 

Let $L_U \subset L$ be a relatively open subset containing $L_V$ in its closure and put $\cF=TL$ on $L_U$. Then for any extension of $\cF$ to a Lagrangian distribution defined on an open neighborhood $L_U \subset U_L \subset K$, we will have $\cH \pitchfork \cF$ on $U_L \cap V$ by the collaring property, possibly after shrinking the $V_x$, thus replacing $V$ with a smaller open set, and then possibly replacing $U_L$ with a smaller neighborhood of $L_U$. So we may set $U= \bigcup_L U_L$ and then the pair $(\cF,\cH)$ is a (dual) $\cI$-structure subordinate to the cover $\cK=U \cup V$. \end{proof}

\begin{figure}[h]
\includegraphics[scale=0.15]{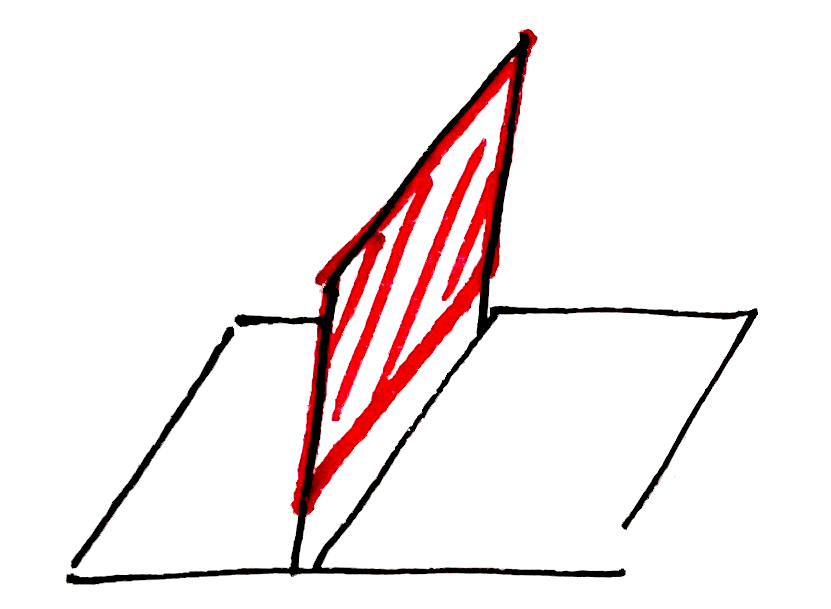}
\caption{The open set $U_L$ for $L$ the smooth piece which is not the zero section in the arboreal singularity $A_2$.}
\label{fig:A_2-U_1}
\end{figure}

\begin{figure}[h]
\includegraphics[scale=0.15]{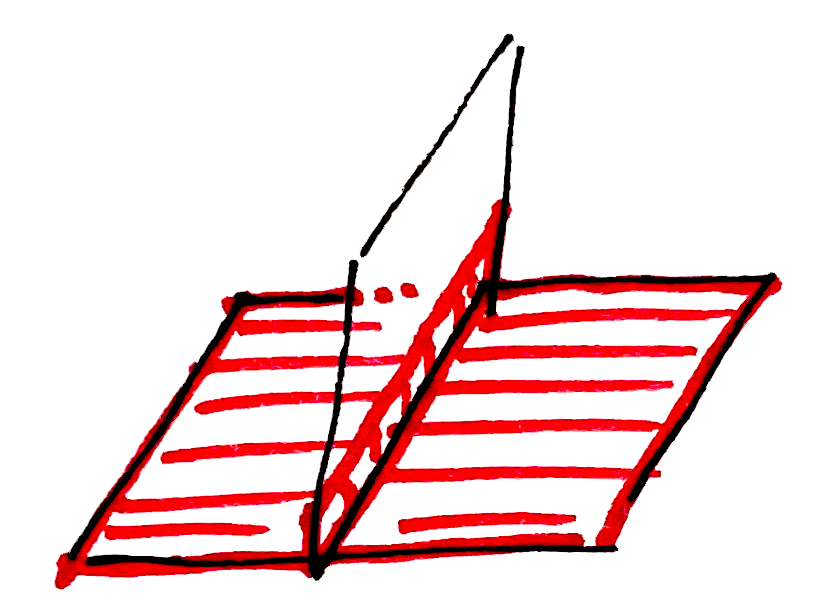}
\caption{The open set $U_L$ for $L$ the zero section of the arboreal singularity $A_2$.}
\label{fig:A_2-U_2}
\end{figure}

\begin{figure}[h]
\includegraphics[scale=0.15]{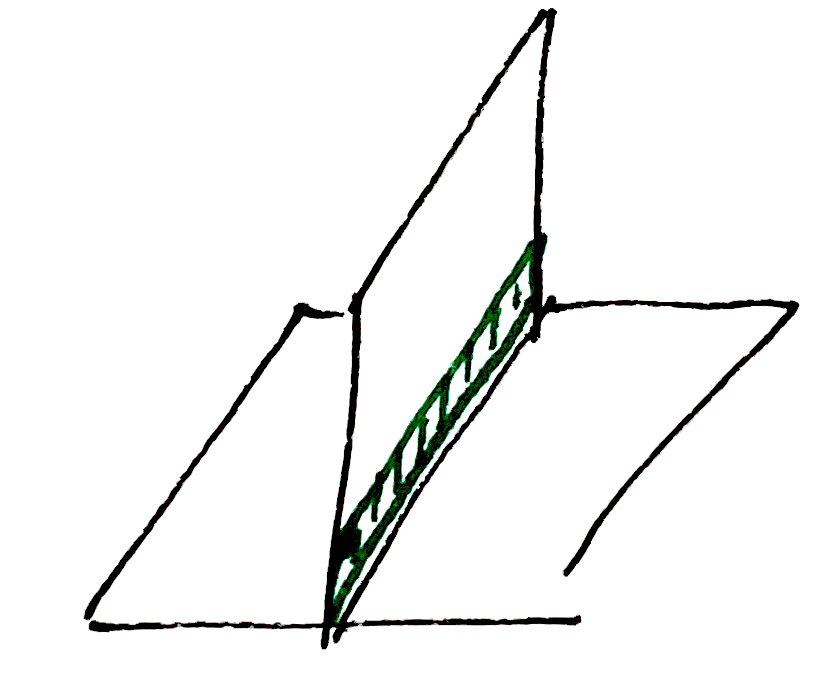}
\caption{The open set $V$ for the arboreal singularity $A_2$.}
\label{fig:A_2-V}
\end{figure}

\begin{figure}[h]
\includegraphics[scale=0.2]{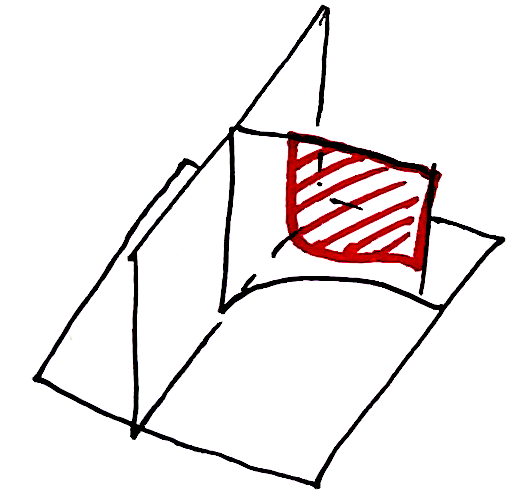}
\caption{The open set $U_L$ for $L$ the smooth piece of the arboreal singularity $A_3$ which is a quarter-space.}
\label{fig:A_3-U_3}
\end{figure}

\begin{figure}[h]
\includegraphics[scale=0.2]{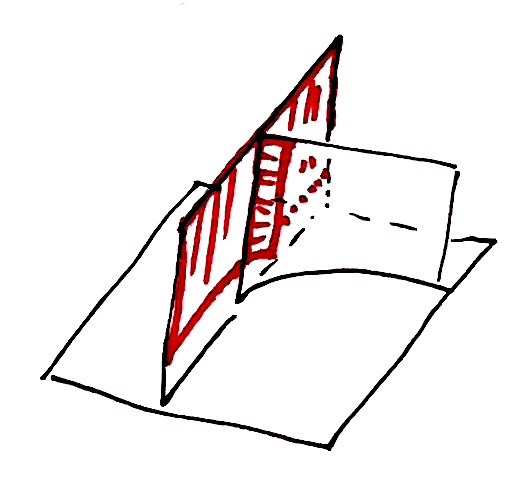}
\caption{The open set $U_L$ for $L$ the smooth piece of the arboreal singularity $A_3$ which is a half-space}
\label{fig:A_3-U_1}
\end{figure}

\begin{figure}[h]
\includegraphics[scale=0.2]{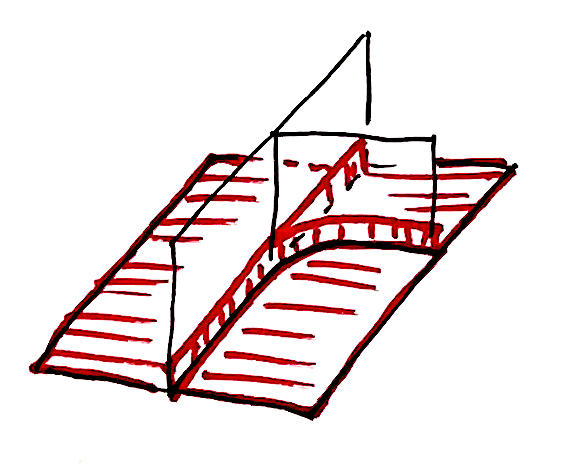}
\caption{The open set $U_L$ for $L$ the zero section of the arboreal singularity $A_3$.}
\label{fig:A_3-U_2}
\end{figure}

\begin{figure}[h]
\includegraphics[scale=0.2]{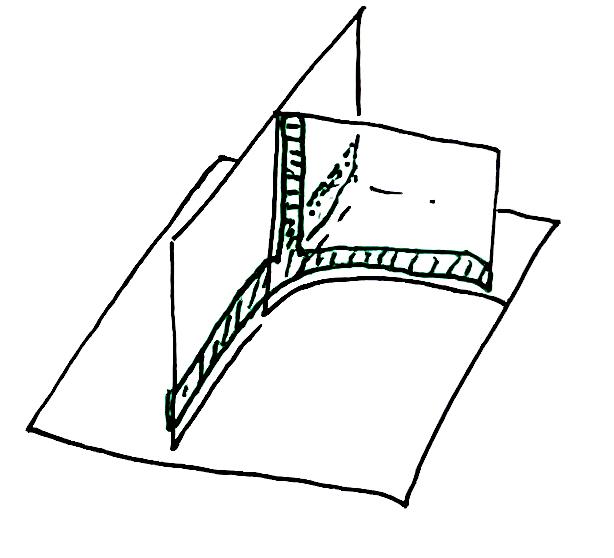}
\caption{The open set $V$ for the arboreal singularity $A_3$.}
\label{fig:A_3-V}
\end{figure}

\begin{remark} See Figures \ref{fig:A_2-U_1}, \ref{fig:A_2-U_2} and \ref{fig:A_2-V} for an illustration of the open cover $U \cup V$ near an $A_2$ arboreal singularity and see Figures \ref{fig:A_3-U_3}, \ref{fig:A_3-U_2}, \ref{fig:A_3-U_1} and \ref{fig:A_3-V} for an illustration of the open cover $U \cup V$ near an $A_3$ arboreal singularity. \end{remark}

\subsection{Cohomological obstructions to existence of arboreal skeleta}

Recall that Theorem \ref{thm:obs-arb} states that if $X$ is a Weinstein manifold which admits an arboreal skeleton, then $c_1(X)c_{2k}(X) = c_{2k+1}(X)$ in $H^{4k+2}(X;\bZ[\frac{1}{2}])$ for all $k\geq 2$.
 
 \begin{proof}[Proof of Theorem \ref{thm:obs-arb}] Follows immediately from Theorem \ref{thm:I-str} and Proposition \ref{prop:lift-obs-main}. \end{proof}

 We will illustrate this cohomological obstruction with examples which are complements of smooth divisors in the complex projective space $\bP^n$. We first recall some basic facts.

Let $Y$ be a smooth projective variety over $\bC$ and $D \subset Y$ be an effective ample divisor cut out by a section $s \in H^0(Y, \cO(D))$. After choosing a metric $\| \cdot \|$ such that the associated curvature $\omega$ is a positive $(1,1)$-form, we obtain a Stein structure on the affine variety $X=Y\setminus D$ by setting the potential $\phi = - \log \| s \|$, see for example \cite{Mc}. In particular, we get a Weinstein structure on $X$, unique up to Weinstein homotopy.

We restrict our attention to the case where $Y=\bP^n$ and $D \subset \bP^n$ is a smooth divisor of degree $d \geq 1$ and consider the Weinstein manifold $X_{n,d} = \bP^n \setminus D$.

We first recall the computation of the cohomology of $X_{n,d}$ with coefficients in a given commutative ring $R$ away from the middle dimension. Note that for $i > n$ we have that $H_i(X_{n,d};R)=0$ since $X_{n,d}$ has the homotopy type of a CW complex of dimension $\leq n$. From this observation alone one may deduce a version of the classical Lefschetz hyperplane theorem, namely that the inclusion $D \subset \bP^n$ induces isomorphisms $H^i(\bP^n;R) \to H^i(D;R)$ when $i<n-1$ for any coefficients $R$.
To see this, simply note that the long exact sequence for the pair $(\bP^n,D)$ reads
$$ \cdots \to H^i(\bP^n ,D;R) \to H^i(\bP^n;R) \to H^i(D;R) \to H^{i+1}(\bP^n,D;R)  \to \cdots$$
By Lefschetz duality $H^i(\bP^n,D;R) \simeq H_{2n-i}(X_{n,d};R)=0$ when $i<n$, from which the desired conclusion follows. 

The well-known computation of the cohomology of $X_{n,d}$ can then be obtained. We record it as a proposition for future reference. 

\begin{proposition}\label{prop:coh-comp}
Let $R \xrightarrow{d} R$ be the multiplication by $d$ map. For $i<n$ odd we have $H^i(X_{n,d};R) \simeq \ker( R \xrightarrow{d} R)$, while for $0<i<n$ even we have $H^i(X_{n,d};R) \simeq \text{coker}( R \xrightarrow{d} R)$, i.e., $R/dR$, and the restriction map $H^i(\bP^n ; R) \to H^i(X_{n,d};R)$ is the surjection $R \to R/dR$.
\end{proposition}

\begin{proof}
Let $[H] \in H^2(\bP^n;\bZ)$ be the hyperplane class, so $H^*(\bP^n;R) \simeq R[H]/[H]^{n+1}$. Consider the long exact sequence for the pair $(\bP^n,U)$  where $U=\bP^n \setminus V$ for $V$ a tubular neighborhood of $D$. This reads
$$ \cdots \to H^i(\bP^n,U;R) \to H^i(\bP^n;R) \to H^i(U;R) \to H^{i+1}(\bP^n,U;R) \to H^{i+1}(\bP^n;R) \to  \cdots $$
By excision $H^i(\bP^n ,U; R) \simeq H^i(V, \partial V;R)$, by the Thom isomorphism $H^i(V, \partial V;R) \simeq H^{i-2}(D;R)$, and by the Lefschetz theorem $H^{i-2}(D;R) \simeq H^{i-2}(\bP^n;R)$. Further, for $i<n$, under this identification the map $H^{i-2}(\bP^n;R) \to H^i(\bP^n;R)$ is given by cup product with the Euler class of the normal bundle to $D$, namely $d[H]$. 

In short, when $0<i<n$ even the exact sequence takes the form
$$ \cdots \to R \xrightarrow{d} R \to H^i(X_{n,d};R) \to 0 \to  \cdots, $$
while for $i<n$ odd the exact sequence takes the form
$$ \cdots 0 \to H^i(X_{n,d};R) \to R \xrightarrow{d} R \to \cdots.  $$
The conclusion follows.
\end{proof}
 
\begin{remark}
In particular, $H^i(X_{n,d};\bZ) \simeq \bZ/d$ for $0<i<n$ even and $H^i(X_{n,d};\bZ)=0$ for $i<n$ odd. The middle dimensional cohomology group $H^n(X_{n,d};\bZ)$ can be quite complicated, but this will not be relevant to us.
\end{remark}

We also recall that for $2i<n$ the $i$-th Chern class $c_i(X_{n,d}) \in H^{2i}(X_{n,d};R) \simeq R/dR$ is the image of ${n+1 \choose i}$ under the ring homomorphism $\bZ \to R/dR$. Indeed, this follows from the above discussion and from the fact that 
\[ c(\bP^n) = \sum_{i=0}^{n} {n+1 \choose i} [H]^i \in H^*(\bP^n; \bZ). \]

 \begin{example}
If $D \subset \bP^n$ is an anti-canonical divisor, i.e., of degree $n+1$, then $c_1(X_{n, n+1})=0$ in $H^2(X;\bZ)$ and so in particular $X_{n, n+1}$ is gradable.
 \end{example}

We can now deduce the following result from Theorem \ref{thm:obs-arb}.

 \begin{cor}\label{cor: basic}
Let $X_{n,d}$ be the complement of a smooth divisor $D$ of degree $d$ in $\bP^n$ where $n > 6$ and $d$ is odd. If the integer $\frac{1}{3}n(n+1)(n+2)$ is not divisible by $d$, then $X_{n,d}$ does not admit an arboreal skeleton. 
\end{cor}

\begin{proof}
By Proposition \ref{prop:coh-comp}, for $i<n$ even have $H^i(X_{n,d};\bZ[\frac{1}{2}]) \simeq \text{coker}( \bZ[\frac{1}{2}] \xrightarrow{d} \bZ[\frac{1}{2}] ) \simeq \bZ/d$ since $d$ is odd. Further, for $2j<n$ the Chern class $c_j(X) \in H^{2j}(X_{n,d};\bZ[\frac{1}{2}])$ corresponds to the residue class ${n + 1 \choose j}\, \mod \, \, d $ in $\bZ/d$. 

Since $n > 6$ by assumption, if $X_{n,d}$ were arborealizable one would have the equality
$$ {n+1 \choose 1}{n+1 \choose 2} - {n+1 \choose 3} \equiv 0 \,  \mod \, \, d $$
which was to be proved.
\end{proof}

\begin{example}
In particular, for fixed $n>6$ it follows that $X_{n,d}$ is never arborealizable when $d$ is a sufficiently large odd integer. 
\end{example}

\begin{example}
The smallest example to which this obstruction applies is $n=7$, $d=3$.
\end{example}

\begin{example}
The integer $\frac{1}{3}n(n+1)(n+2)$ is divisible by $n+1$ if and only if $n$ is congruent to 0 or 1 mod 3. Hence for any  $n \geq 8$ congruent to 2  mod 6 the complement $X$ of a smooth anti-canonical divisor in $\bP^{n}$ is a gradable Weinstein manifold which does not admit an arboreal skeleton. The smallest example to which this obstruction applies is $n=8$.
\end{example}

\section{The Maslov obstruction}\label{sec: Maslov}

\subsection{Maslov data}

A stably presentable symmetric monoidal $\infty$-category $\cC$ \cite{Lutopos}, \cite{Lualg} has an associated group $\mathcal{P}ic(\cC)$ of isomorphism classes of invertible objects called its {\em Picard group}. Classical examples include the category of modules over a commutative ring and the category of quasi-cohorent modules over the structure sheaf of a scheme. We have $\cP ic(\cC)=\pi_0 (\Pic(\cC))$ for an infinite loop space $\Pic(\cC)$ (the zeroth space of a connective spectrum called the Picard spectrum of $\cC$). When $\cC$ is the category of modules over a commutative ring or commutative ring spectrum $R$ we write $\Pic(R)$. There holds $\Omega \Pic(R) \simeq \text{GL}_1 (R)$ and so we have a fibration $$B \text{GL}_1(R)\to \Pic(R) \to \pi_0 (\Pic(R)).$$ 
In particular there holds $\pi_1 (\Pic(R)) \simeq \pi_0(\text{GL}_1(R))$ and $\pi_k(\Pic(R)) \simeq \pi_{k-1}(R)$ for $k \geq 2$. For example when $R=\bZ$ we have a fibration $B \bZ/2 \to \Pic(\bZ) \to \bZ$ and so $\pi_0( \Pic(\bZ))=\bZ$, $\pi_1( \Pic(\bZ)) =\bZ/2$ and $\pi_k( \Pic(\bZ)) = 0$, $k \geq 2$. 

When $R$ is the sphere spectrum $\bS$ there holds $\pi_0( \Pic(\bS)) = \bZ$ \cite{HMS} and we have a fibration $$B \text{GL}_1(\bS) \to \Pic(\bS) \to \bZ. $$ In particular we have $\pi_k( \Pic(\bS)) \simeq \pi_{k-1}^s$ for $k\geq2$, the stable homotopy groups of spheres, and $\pi_1(\Pic (\bS)) =\bZ/2$. This case is universal in that for any other $\cC$ we have a natural infinite loop map $\Pic(\bS) \to \Pic(\cC)$. 

The Maslov obstruction of a Weinstein manifold $X$ is the composition
$$X  \to BU \to B(U/O) \to  B^2\Pic(\bS) $$
where the map $X \to BU$ classifies $TX$ and the map $B(U/O) \to  B^2\Pic(\bS)$ is a universal map given by microlocal sheaf theory \cite{NS},  which according to \cite{J1,J2} is a three-fold delooping of the real $J$-homomorphism 
$$ J: O \to \text{GL}_1(\bS). $$
More precisely, the fibrations $BO \to \Omega(U/O) \to \bZ$ and $B \GL_1(\bS) \to \Pic(\bS) \to \bZ$ deloop to fit into a homotopy commutative diagram in which the columns are fibrations
\[ \xymatrix{  B^3O \ar[d]\ar[r]^{B^3J} & B^3 \GL_1( \bS ) \ar[d] \\ B(U/O) \ar[r]  \ar[d] & B^2\Pic(\bS) \ar[d] \\ B^2\bZ \ar[r]^{=} & B^2 \bZ } \]

See \cite{J1,J2} for a concrete model of the once delooped real $J$-homomorphism  $BJ: \bZ \times BO \to \Pic(\bS)$, which classifies the spherical fibration over $\bZ \times BO$ given by the fiberwise one-point compactification of the tautological bundle.

For any other $\cC$, we refer to the map $X \to B^2\Pic(\cC)$ given by the composition of $X \to B^2\Pic(\bS)$ and the natural map $B^2\Pic(\bS)  \to B^2\Pic(\cC)$ as the $\cC$-Maslov obstruction. 

\begin{definition} When the Maslov obstruction (resp. $\cC$-Maslov obstruction) of a Weinstein manifold $X$ vanishes, we say that $X$ admits Maslov data (resp. $\cC$-Maslov data), and a choice of nullhomotopy is called a choice of Maslov data (resp. $\cC$-Maslov data).  \end{definition}

Given a choice of $\cC$-Maslov data, Nadler and Shende construct in \cite{NS} a microlocal category for $X$ over $\cC$, invariant under Weinstein homotopy.

\begin{example}\label{ex:gradable} Let us explain why the condition $2c_1(X)=0$ in $H^2(X;\bZ)$ is an obstruction to the existence of $\bZ$-Maslov data (and hence also to the existence of $\bS$-Maslov data). Delooping the fibration $B(\bZ/2) \to \Pic(\bZ) \to \bZ$ we obtain a fibration $B^{3} (\bZ /2 )\to B^2\Pic(\bZ) \to B^2\bZ=BU(1)$ where the composition $BU \to B^2\Pic(\bZ) \to B^2\bZ$ is identified as the map $BU \to BU(1)$ given by $B(\text{det}^2)$, see \cite{G, CKNS}. So the vanishing of the $\bZ$-Maslov obstruction $X \to B^2\Pic(\bZ)$ forces $2c_1(X)=0$ in $H^2(X;\bZ)$, which is of course the familiar obstruction to the construction of $\bZ$-gradings in the Fukaya category.  \end{example}

\begin{example} If $X$ admits a polarization, i.e., if $X \to BU$ lifts to a map $X \to BO$, then the composition $X \to BO \to BU \to B(U/O)$ is nullhomotopic. So if $X$ admits a polarization, then the $\cC$-Maslov obstruction of $X$ vanishes  for every $\cC$. Further, given a choice of lift $X \to BO$ we get a canonical choice of  nullhomotopy for the composition $X \to B(U/O)$ and hence a canonical choice of $\cC$-Maslov data over every $\cC$. \end{example}

\subsection{Cohomological obstructions to the existence of Maslov data} 

In this section we establish secondary obstructions to the existence of Maslov data on a Weinstein manifold $X$, beyond the gradability obstruction $2c_1(X) =0$ in $H^2(X;\bZ)$ explained in Example \ref{ex:gradable}.

\begin{theorem}\label{thm: char class}
For each odd prime $p$ there exists a characteristic class for complex vector bundles $m_p \in H^{2p}(BU;\bZ/p)$ satisfying the following properties:
\begin{enumerate}
\item If $X$ is a Weinstein manifold which admits Maslov data, then $m_p(X)=0$.
\item $m_p$ is a polynomial in the Chern classes $c_1,c_2,\ldots,c_p$.
\item The $\bZ/p$ coefficient of $c_p$ in $m_p$ is nonzero.
\end{enumerate}
\end{theorem}

Here $m_p(X) \in H^{2p}(X;\bZ)$ denotes the pullback of $m_p$ under the classifying map $X \to BU$ of $TX$. An immediate consequence is Theorem \ref{thm:obs-mas}, which we recall states that if $X$ is a Weinstein manifold which admits Maslov data, and if for some odd prime $p$ we have that $c_1(X)=c_2(X)= \cdots = c_{p-1}(X)=0$ in $H^*(X;\bZ/p)$, then we must also have $c_p(X)=0$ in $H^*(X;\bZ/p)$.

\begin{proof}[Proof of Theorem \ref{thm:obs-mas}]
If $c_i(X)=0$ in $H^{2i}(X;\bZ/p)$ for $1\leq i\leq p-1$, then by properties (ii) and (iii) of Theorem \ref{thm: char class} it follows that $m_p(X)$ is a nonzero multiple of $c_p(X)$. If $X$ admits Maslov data, then by property (i) of Theorem \ref{thm: char class} we have $m_p(X)=0$, and so we conclude that $c_p(X)=0$ as desired.\end{proof}

\begin{proof}[Proof of Theorem \ref{thm: char class}]
We will use the following  fact about the stable homotopy groups of spheres $\pi_k^s$ due to work of Adams and Quillen: for each odd prime $p$, we have that $\pi_k^s \otimes \bZ/p = 0$ for all $ 0 < k < 2p-3$. Furthermore, $\pi_{2p-3}^s \otimes \bZ/p \simeq \bZ/p$ and this single $\bZ/p$ factor is in the image of the real J-homomorphism $\pi_{2p-3}(O) \to \pi_{2p-3}^s$ (see \cite{Ra} for a reference). Since we have $\pi_{*} (\GL_1(\bS))= \pi_{*}^s$, we can conclude from the mod $p$ Hurewicz theorem \cite{Se} that
\[ H^{2p}(B^3 \GL_1(\bS); \bZ/p) \simeq \pi_{2p}(B^3 \GL_1(\bS)) \otimes \bZ/p  \simeq  \bZ/p.\] 
From the Serre spectral sequence associated to the fibration $B^3 \GL_1(\bS) \to B^2 \Pic(\bS) \to B^2 \bZ$ we further conclude that 
\[ H^{2p}(B^2 \Pic (\bS); \bZ/p) \simeq H^{2p}(B^3 \GL_1(\bS); \bZ/p) \simeq  \bZ/p \] 
with the isomorphism induced by the map $B^3 \GL_1(\bS) \to B^2 \Pic(\bS)$.

Now let $\alpha$ denote a generator of $H^{2p}(B^2 \Pic (\bS); \bZ/p)$, and let $m_p$ denote the pullback of $\alpha$ to $H^{2p} (BU; \bZ/p)$, so $m_p$ is some polynomial in the Chern classes with total degree $2p$, and for a Weinstein manifold $X$ the pullback of $m_p$ to $H^{2p}(X;\bZ/p)$ is a class $m_p(X)$ which obstructs the triviality of the Maslov obstruction $X \to B^2\Pic(\bS)$.  

Let $\lambda \in \bZ/p$ denote the coefficient of $c_{p}$ in $m_p$. To show that $\lambda \neq 0$, consider a complex vector bundle $\cE$ over $S^{2p}$ classified by a generator of $\pi_{2p}(BU)$. We claim that the full composite \[f: S^{2p} \to BU \to B(U/O) \to B^2 \Pic(\bS)\] must generate $\pi_{2p}( B^2 \Pic(\bS)) \otimes \bZ/p$. First, since $B(U/O) \to B^2\Pic(\bS)$ is a delooping of the J-homomorphism, by the discussion above we see that the map $\pi_{2p}(B(U/O))  \otimes \bZ/p \to \pi_{2p} (B^2 \Pic(\bS)) \otimes \bZ/p$ surjects. Thus 
to establish the claim it suffices to know that \[\pi_{2p} (BU ) \otimes \bZ/p \to \pi_{2p}( B(U/O)) \otimes \bZ/p\] is surjective. This is straightforward: the portion
$$\pi_{2p-1} (U) \to \pi_{2p-1} (U/O) \to \pi_{2p-2} (O)\to 0$$
of the long exact sequence in homotopy groups associated to the fibration $O \to U \to U/O$ can be directly computed from Bott periodicity.  Since $2p-1$ is odd and cannot take the residue classes 3 or 7 mod 8, the sequence is either
$$ \bZ \xrightarrow{2} \bZ \to \bZ/2  \to 0 \quad \text{or} \quad \bZ \xrightarrow{\simeq} \bZ \to 0 \to 0;$$
the desired conclusion follows in either case. Hence we have shown that the pushforward
\[f_* : \pi_{2p} (S^{2p}) \otimes \bZ/p \to \pi_{2p} (B^2 \Pic(\bS)) \otimes \bZ/p\]
is an isomorphism. Again by the mod $p$ Hurewicz theorem, this means the pullback on cohomology
\[ f^*: H^{2p}(B^2 \Pic (\bS); \bZ/p) \to H^{2p}(S^{2p}; \bZ/p)\]
 is also an isomorphism (to be more precise: one argues with mod $p$ Hurewicz that this is true for the lift $\tilde{f}: S^{2p} \to B^3 \GL_1(\bS)$, which exists since $\pi_{2p}(B^2 \bZ)=0$, and this implies the result for $f$). Thus the pullback of $m_p$ is nonzero in $H^{2p}(S^{2p}; \bZ/p)$. Of course, for dimension reasons the lower Chern classes are all zero in $S^{2p}$, so it has to be some nonzero multiple of $c_p$, which shows that $\lambda \neq 0$ as desired. \end{proof}

\begin{remark} We do not further determine the coefficients of the characteristic class $m_p$ here and will content ourselves with the elementary consequence Theorem \ref{thm:obs-mas} for our sample applications.
\end{remark}

\begin{corollary}\label{cor: Fermat}
The complement $X$ of a smooth anti-canonical divisor in $\bP^{kp-1}$ does not admit Maslov data for any prime $p>2$ and any $k \geq 3$ not divisible by $p$. 
\end{corollary}

\begin{proof}
By Proposition \ref{prop:coh-comp} the cohomology $H^j(X;\bZ/p)$ is isomorphic to $\bZ/p$ for $j<kp-1$ and when $2i<kp-1$ the Chern class $c_i(X)$ is the residue class ${kp \choose i}$ in $\bZ/p$. In particular this applies for $i=1,\ldots,p$. It is straightforward to check that for $i<p$ the integer ${kp \choose i}$ is divisible by $p$ while ${kp \choose p}$ is not. The conclusion then follows from Theorem \ref{thm:obs-mas}. \end{proof}

\begin{example}  The smallest example $n,d$ to which this obstruction applies is $\bP^{11} \setminus D_{12}$, with $k=4$, $p=3$.  \end{example}

\begin{question}
Complements of smooth anti-canonical divisors  in $\bP^n$ are gradable ($2c_1=0$), have integer Fukaya categories and satisfy various formulations of homological mirror symmetry over the integers. However, from Corollary \ref{cor: Fermat} we see that some of them {\em shouldn't} have spectral Fukaya categories.  Do their integer Fukaya categories have any algebraic property that reflects this?
\end{question}

\subsection{Vanishing of the Maslov obstruction does not imply arborealizability}\label{sec:crazy-ex}

In this final section we show that existence of Maslov data does not imply existence of an arboreal skeleton (or of a polarization). 

\begin{proposition}\label{prop: crazy-ex}
There exists a Weinstein manifold $X$ such that:
\begin{enumerate}
\item $X$ is homotopy equivalent to $S^6$.
\item $c_3(X) \neq 0$ in $H^6(X;\bZ) \simeq \bZ$.
\item The image of a generator of $\pi_6(X) \simeq \bZ$ under the homomorphism $\pi_6(X) \to \pi_6(BU)$ induced by the classifying map $X \to BU$ of $TX$ is a multiple of $24$ in $\pi_6(BU) \simeq \bZ$.
\end{enumerate}
\end{proposition}

Before giving the proof of Proposition \ref{prop: crazy-ex} we deduce Corollary \ref{cor: ex Maslov}.
 
\begin{proof}[Proof of Corollary \ref{cor: ex Maslov}] Let $X$ be as in Proposition \ref{prop: crazy-ex}. Such an $X$ cannot be arborealizable by Theorem \ref{thm:obs-arb}. To see that the Maslov obstruction vanishes, we only need to check that it induces the zero homomorphism on $\pi_6$. We have seen that the map $\pi_6(B (U/O)) \to \pi_6( B^2\Pic(\bZ))$ can be identified with the J-homomorphism $\pi_3 (O) \to \pi_3^s$, which is the surjection $\bZ \to \bZ/24$. Since by hypothesis the map $X \to BU$ maps a generator of $\pi_6 (X) \simeq \bZ$ to an element of $\pi_6( BU) \simeq \bZ$ which is a multiple of $24$, it follows that $\pi_6 (X) \to \pi_6( B^2 \Pic(\bS))$ is the zero homomorphism, and hence that $X \to B^2 \Pic(\bS)$ is nullhomotopic as desired. 
\end{proof}

\begin{remark}
 Note that $H^6(X;\bZ/3) \simeq \bZ/3$ and the Chern class $c_3(X) \in H^6(X;\bZ)$ is a multiple of $24$, hence vanishes in $H^6(X;\bZ/3)$ as required by Theorem \ref{thm:obs-mas}.
 \end{remark}

\begin{proof}[Proof of Proposition \ref{prop: crazy-ex}.] Let $p:\cV \to S^6$ denote the complex vector bundle over $S^6$ given by $TS^6$ with the homotopically unique almost complex structure on $S^6$. Note that $c_3(\cV)=e(TS^6)$ is nonzero in $H^6(S^6;\bZ) \simeq \bZ$, indeed $c_3(\cV) \cdot [S^6]=\chi(S^6)=2$ for $[S^6] \in H_6(S^6;\bZ)$ the fundamental class of $S^6$. 

Next, let $\pi: \cW \to S^6$ be another complex vector bundle, and let $X$ denote its total space. We may (non-canonically) split $TX \simeq \pi^*\cW \oplus \pi^*TS^6$ and thereby equip $X$ with an almost complex structure $J$ by taking the direct sum of the complex structures on $\cW$ and $\cV$. 

By the Eliashberg h-principle on existence of Stein structures \cite{E}, if the complex rank of $\cW$ is at least $3$ there exists a Weinstein structure on $X$ compatible with $J$. Indeed, there exists a Stein structure on $X$ whose complex structure $J'$ is homotopic to $J$ in the space of almost complex structures on $X$. 

For $X$ equipped with this Weinstein structure we therefore have $c_3(X)=\pi^*c_3(\cW)+\pi^*c_3(\cV)$. If we take $\cW=  \cV^{\oplus k}$, the direct sum of $k \geq 1$ copies of $\cV$, we will certainly have $c_3(X) \neq 0$ in $H^6(X;\bZ) \simeq \bZ$. Indeed, there holds 
$$c_3(X) \cdot \iota_* [S^6] =  (k+1) \pi^* c_3(\cV) \cdot \iota_*[S^6] =  (k+1) c_3(\cV) (\pi_* \iota_*)  [S^6] = (k+1)c_3(\cV) \cdot [S^6] = 2(k+1) \neq 0$$ 
for $[S^6] \in H_6(S^6;\bZ)$ the fundamental class of $S^6$ and $\iota: S^6 \hookrightarrow X$ the inclusion of the zero section.

Let $\alpha$ (resp. $\beta$) in $\pi_6(BU)$ denote the image of the generator $[S^6 \xrightarrow{\text{id}} S^6] \in \pi_6( S^6) \simeq \bZ$ under the homomorphism $\pi_6( S^6) \to \pi_6( BU)$ induced by the map $S^6 \to BU$ which classifies $\cV$ (resp. $\cW$). Then the classifying map $X \to BU$ of $TX$ maps the generator $[\iota:S^6 \hookrightarrow X] \in \pi_6(X) \simeq \bZ$  to the sum of $\alpha$ and $\beta$ in $\pi_6(BU)$. Further, if $\cW=\cV^{ \oplus k }$, then $\beta= k \alpha$ and the classifying map $X \to BU$ of $TX$ maps the generator $[\iota] \in \pi_6(X)$ to $(k+1) \alpha$. It therefore suffices to pick $k=23$. \end{proof}

We conclude with some final remarks.

\begin{remark}\label{rem: final} The Weinstein manifold $X$ constructed above has real dimension 144 and is certainly very subcritical. However, by a theorem of Cieliebak \cite{Cie02} it is deformation equivalent to a Weinstein manifold of the form $Y \times \bC^{66}$ where $Y$ is a 12-dimensional Weinstein manifold which is not subcritical. The destabilized Weinstein manifold $Y$ also satisfies the properties of Proposition \ref{prop: crazy-ex}.  Alternatively, since $\pi_k(U(m)) \to \pi_k(U(m+1))$ is an isomorphism for $k<2m$ and surjective for $k=2m$, we see that the classifying map $S^6 \to BU(69)$ of the vector bundle $\cV^{\oplus 23}$ factors as a composition $S^6 \to BU(3) \to BU(69)$, i.e., $\cV^{\oplus 23}$ is isomorphic as a complex vector bundle to a stabilization $\cY \oplus \underline{\bC}^{66}$ where $\cY \to S^6$ is a complex vector bundle of complex rank 3. The total space of $\cY$ is an almost complex manifold of real dimension 12, which may be endowed with a Weinstein structure by the h-principle \cite{E} to produce the desired critical Weinstein manifold $Y$ satisfying the properties of Proposition \ref{prop: crazy-ex}.  \end{remark}

\begin{remark} Note that in any event our $Y$ is produced by the h-principle and therefore a priori is flexible. If one wanted an example of a non-flexible Weinstein manifold as in Corollary \ref{cor: ex Maslov} one could take the boundary connect sum $Z$ of $Y$ and $T^*S^6$, which has a nontrivial Fukaya category (a subcritical handle attachment does not change the Fukaya category). To see that $Z$ satisfies the conditions of Corollary \ref{cor: ex Maslov}, note that the homotopy type of $Z$ is that of a wedge of two 6-spheres and so the map $Z \to BU$ may be identified with the wedge of the two maps $S^6 \to BU$ obtained by restricting $TZ$ to the zero sections of $T^*S^6$ and $Y$. Since the Maslov obstructions of $Y$ and $T^*S^6$ are trivial ($T^*S^6$  is polarizable), we conclude that the Maslov obstruction of $Z$ is trivial, i.e. that $Z$ admits Maslov data.  However, we have $c_3(Z) \neq 0$ in $H^6(Z;\bZ[\frac{1}{2} ] ) \simeq \bZ[\frac{1}{2} ]  \oplus \bZ[\frac{1}{2} ] $ since $c_3(Z)$ still pairs to $48$ with the fundamental class of the zero section of $Y$, so $Z$ does not admit an arboreal skeleton by Theorem \ref{thm:obs-arb}.  \end{remark}

\begin{question}
What does the skeleton of $X$ (or $Y$, or $Z$) look like? In particular, what are the simplest Lagrangian singularities beyond arboreal singularities which are needed?
\end{question}

\printbibliography

 \end{document}